# A PHASE TRANSITION FOR COMPETITION INTERFACES

BY PABLO A. FERRARI, JAMES B. MARTIN AND LEANDRO P. R. PIMENTEL

*Universidade de São Paulo, University of Oxford and Delft University of Technology*

We study the *competition interface* between two growing clusters in a growth model associated to last-passage percolation. When the initial unoccupied set is approximately a cone, we show that this interface has an asymptotic direction with probability 1. The behavior of this direction depends on the angle $\theta$ of the cone: for $\theta \geq 180°$, the direction is deterministic, while for $\theta < 180°$, it is random, and its distribution can be given explicitly in certain cases. We also obtain partial results on the fluctuations of the interface around its asymptotic direction. The evolution of the competition interface in the growth model can be mapped onto the path of a *second-class particle* in the totally asymmetric simple exclusion process; from the existence of the limiting direction for the interface, we obtain a new and rather natural proof of the strong law of large numbers (with perhaps a random limit) for the position of the second-class particle at large times.

**1. Introduction.** The behavior of the boundary (or "growth interface") of a randomly growing cluster has been much investigated. In many models, the growing region, after linear rescaling, is seen to converge to a deterministic asymptotic shape, and the fluctuations of the growth interface follow either Gaussian or Kardar–Parisi–Zhang scales. See, for example, [19] for a wide-ranging review.

A less well-studied phenomenon is the "competition interface" between two clusters growing in the same space. Derrida and Dickman [4] describe simulations of a first-passage percolation model (or Eden model) in which two clusters grow into a vacant sector of the plane with angle $\theta$. They obtain values for the roughening exponents of the competition interface: the fluctuations of this interface at distance $r$ are of order $r^\chi$, where $\chi = 1/3$ for $\theta < 180°$, $\chi = 2/3$ for $\theta = 180°$ and $\chi = 1$ for $\theta > 180°$. However, they note that $\chi = 1$ does not describe the true roughness of the interface in the last case, but instead indicates a *random direction*; if instead $\chi$ is defined in terms of the fluctuations of the interface about its (maybe random) asymptotic direction, one should expect to see $\chi = 2/3$.

We investigate analogous questions for a related growth model which is associated to *directed last-passage percolation* in the plane. The exact solvability of this model makes it possible to obtain rigorous results about the existence and distribution of an asymptotic direction of the interface, and about its roughness. A further important motivation is the relation between the competition interface in this model and the behavior of a second-class particle in the totally asymmetric simple exclusion process (TASEP), which we describe below.

The growth model can be described as follows (see Section 2 for precise definitions and notation). At a given time $t$ the occupied set $\Gamma_t$ is a decreasing subset of $\mathbb{Z}^2$ (i.e., if $z$ is occupied, then so is any point below and to the left of $z$). If $z$ is unoccupied but both $z - (0,1)$ and $z - (1,0)$ are occupied, then $z$ is added to the occupied set at rate 1. That is, if $G(z)$ is the time at which $z$ joins the occupied set, then the quantities $G(z) - \max\{G(z - (0,1)), G(z - (1,0))\}$ are exponential random variables with rate 1 (and independent for different $z$). Since the occupied set $\Gamma_t$ is a decreasing set, we can equivalently consider the growth interface $\gamma_t$ which is the boundary of $\Gamma_t$, and which is a path in $\mathbb{Z}^2$ taking steps down and to the right.

We consider initial growth interfaces which pass through the points $(-1, 0)$, $(0, 0)$ and $(0, -1)$. Hence $\gamma_0$ consists of one (half-infinite) path in the lower-right quadrant and another in the upper-left quadrant. The competition occurs as follows. Each point in the upper-left part of $\gamma_0$ belongs to cluster 1, and each point in the lower-right part belongs to cluster 2. When a new point $z$ is added to the occupied set, it joins the same cluster as $\tilde{z}$, where $\tilde{z}$ is the argument that maximizes $G(z - (0,1))$ and $G(z - (1,0))$. That is, it joins the same cluster as whichever of its neighbors (below and to the left) was occupied most recently. [The label of the site $(0,0)$ may be left ambiguous, but we stipulate that site $(0,1)$ always joins cluster 1 and site $(1,0)$ always joins cluster 2.] Let $\Gamma_t^1$ and $\Gamma_t^2$ be the points that are occupied at time $t$ by cluster 1 or cluster 2, respectively (so $\Gamma_t = \Gamma_t^1 \cup \Gamma_t^2$). We can further consider the sets $\Gamma_\infty^1$ and $\Gamma_\infty^2$ of points which (eventually) join clusters 1 and 2, respectively. These sets are separated by the *competition interface*, which is a directed path with up-right steps, lying in the positive quadrant (Figure 1).

We will assume that both ends of the initial growth interface $\gamma_0$ have an asymptotic direction. Thus the unoccupied set at time 0 is approximately a cone, with some angle $\theta \in [90°, 270°]$. We will show that under this assumption, the competition interface converges with probability 1 to a limiting direction. There is a phase transition in the behavior as a function of $\theta$. For concave or flat initial sectors ($\theta \geq 180°$), the limiting direction is deterministic (and we calculate it); for a convex sector ($\theta < 180°$) the limiting direction is random. (The phases are opposite to those observed by Derrida and Dickman, since we consider last-passage rather than first-passage growth rules.)

The last-passage percolation representation is as follows: $G(z)$ can be described as the maximal weight of a directed (up-right) path starting from some point in $\gamma_0$ and ending at $z$, where the weights at each site of $\mathbb{Z}_+^2$ are i.i.d. exponential random variables with mean 1. Since the initial interface $\gamma_0$ is the union of two half-lines with asymptotic directions, we are led to develop a sequence of results concerning a *point to half-line percolation* model, in order to prove the existence of the asymptotic direction of the competition interface. In particular, we extend the concept of "$h$-straightness" for the forest of optimal paths, developed for point-to-point percolation models by Newman and co-authors in [22] and later papers.

One of the main motivations is the relation between the competition interface and the second-class particle in the TASEP, developed by Ferrari and Pimentel [11] to treat the case where the initial unoccupied sector is exactly the positive quadrant. The TASEP has state space $\{0,1\}^\mathbb{Z}$. At each site of $\mathbb{Z}$, there is either a particle or a hole. Each particle tries to jump at rate 1 to the right, and a jump succeeds if the site to its right is unoccupied (so that the particle and hole then exchange places). A *second-class particle* arises when the initial configuration is modified at a single site, and the two processes with and without the perturbation are allowed to evolve using the same random mechanism (the so-called *basic coupling*). At later times the two processes still differ at exactly one site, and the position of this discrepancy is called a second-class particle. More concretely, it moves in the TASEP according to the following rule: it jumps to the right at rate 1 when the site on its right is unoccupied, and when the site on its left contains a particle, it exchanges places with this particle at rate 1.

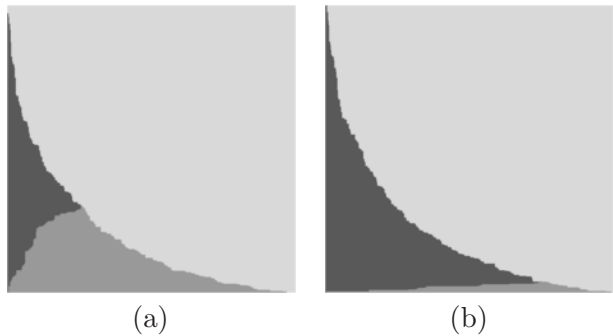

FIG. 1. *Numerical simulations of the competing growth model with $\theta = 90°$, where the unoccupied set at time 0 is precisely the positive quadrant. The pictures show the clusters $\Gamma_t^1$ and $\Gamma_t^2$ at $t = 100$, for two different realizations* (a) *and* (b). *The darker cluster is $\Gamma_t^1$ and the gray one is $\Gamma_t^2$; background is light gray.*

The TASEP and the last-passage percolation growth model can be coupled, using a construction that originates with Rost [25]. Using the construction from [11], the TASEP with a second-class particle can be represented by a TASEP with an extra site, in which the site containing the second-class particle is replaced by a pair of sites containing a hole on the left and a particle on the right. Thus the TASEP with a second-class particle can also be coupled with the growth process, and one obtains in fact that the path of the second-class particle corresponds to the evolution of the competition interface. Under this coupling, initial growth interfaces $\gamma_0$ satisfying the conditions described above correspond to initial TASEP configurations $\eta_0$ with a second-class particle at the origin, and in which the configurations to the left and to the right of the origin have some asymptotic particle densities, say $\lambda$ and $\rho$, respectively. The limit result for the direction of the competition interface can be shown to correspond to a strong law of large numbers for the position of the second-class particle. If $X_t$ is the position of the second-class particle at time $t$, then $X_t/t$ converges with probability 1, to a limit which is deterministic if $\lambda \leq \rho$ (corresponding to the case $\theta \geq 180°$) and random if $\lambda > \rho$. See the discussion after Theorem 3 for references to versions of these laws of large numbers for $X_t/t$ which have previously appeared in the literature.

A particularly interesting case occurs when $\lambda > \rho$ and the initial TASEP configuration is distributed according to product measure with particle density $\lambda$ to the left of the origin and density $\rho$ to the right. In this case the random limit of $X_t/t$ is precisely the uniform distribution on $(1-2\lambda, 1-2\rho)$. This had been obtained as a limit in distribution by Ferrari and Kipnis [8], and the almost sure convergence was proved by Mountford and Guiol [21]. Using the distributional limit one can derive the distribution of the limiting direction of the competition interface in the corresponding case, where the two parts of the initial growth interface $\gamma_0$ are given by two independent random walks going up-left and down-right from the origin. Indeed, the interplay between the two equivalent models of TASEP and growth process is particularly satisfying in this case. Our methods using the competition interface provide a new and rather intuitive proof of the almost sure convergence of $X_t/t$ to some limit; however, the distribution of this limit (and hence also the distribution of the limiting direction of the competition interface) is most naturally derived in the context of the TASEP as was done in [8].

We also provide some partial results about the fluctuations of the competition interfaces. When the initial interface is of random-walk type and the initial growth interface is concave, we show that the roughening exponent of the interface is $\chi = 1/2$. In this case we show that the behavior of the interface on the scale of the fluctuations can effectively be read off from the initial growth interface, and we compute explicitly the covariance matrix of the bidimensional Gaussian vector describing the fluctuations about the asymptotic inclination. For flat initial sectors of random-walk type, we show that $\chi = 2/3$; this is done by showing that the competition interface has the same distribution as an infinite geodesic in the last-passage percolation model, whose roughening exponent can be derived from the results of Balázs, Cator and Seppäläinen [1]. For convex initial sectors we can show that $\chi \leq 3/4$, by an argument involving the bounding of the competition interface between two infinite geodesics with the same asymptotic direction.

We now comment briefly on other recent results concerning related competition growth models. In a first-passage model in which each cluster starts from, say, a single grain, it is already possible that one cluster is surrounded by the other and is unable to grow further. Various results showing that mutual unbounded growth is possible have been obtained in, for example, [3, 12, 13, 15]. In these models on $\mathbb{Z}^d$, the lack of information about the asymptotic shape typically makes it impossible to prove that the competition interface has a limiting direction; in models with rotational symmetry, however, more is possible; see, for example, [23] for results in the case of a model based on a Delaunay triangulation. For the polynuclear growth model, which is closely related to last-passage percolation, results on asymptotic directions and fluctuations for competition interfaces and similar objects are obtained in [2]. Similar phenomena are also observed in [17, 18] for competition interfaces in models which include voter-type dynamics as well as growth, so that the cluster to which an occupied site belongs may not be constant over time.

The rest of the paper is organized as follows. In Section 2 we describe the models precisely, and state the main results. Theorem 1 concerns the existence of a limiting direction for the competition interface and Theorem 2 concerns the distribution of this direction. The corresponding results for the TASEP are given in Theorem 3 (followed by references on related results which already exist in the literature for the TASEP). Our results on fluctuations are given in Theorems 4–6. During Section 2 we also give a "shape theorem" for the percolation model (Proposition 2.1) and explain how the different possible behaviors of the asymptotic shape correspond to the different phases observed for the limiting direction of the competition interface. In Section 3 we prove a sequence of results concerning the point to half-line last-passage percolation model, and use them to prove the laws of large numbers for the competition interface. The proofs of fluctuation results are given in Section 4.

**2. Definitions and results.**

2.1. *Last-passage percolation and competition interfaces.* Let $\mathcal{X} := \{X(z), z \in \mathbb{Z}^2\}$ be a family of i.i.d. exponential random variables with mean 1. Let $\mathbb{P}$ be the probability induced by these variables. For $z \in \mathbb{Z}^2$, we write $z = (z(1), z(2))$ and $|z| = |z(1)| + |z(2)|$. A *directed path* in $\mathbb{Z}^2$ is a (finite or infinite) path $z_1, z_2, \ldots$ such that $z_{i+1} - z_i \in \{(0,1), (1,0)\}$ for each $i = 1, 2, \ldots$. For $z \leq z'$, let $\Pi(z, z')$ be the set of directed paths starting at $z$ and ending at $z'$. The *last-passage time* from $z$ to $z'$ is defined by

$$(2.1) \qquad T(z \to z') = \max_{\pi \in \Pi(z,z')} \left\{ \sum_{y \in \pi} X(y) \right\}.$$

The path $\pi$ realizing the maximum is called the *geodesic* connecting $z$ to $z'$. A semi-infinite path $z_0, z_1, \ldots$ is called a semi-infinite geodesic if for all $k, \ell$, the subpath $z_k, \ldots, z_\ell$ is the geodesic connecting $z_k$ and $z_\ell$. Geodesics satisfy the backward recurrence property

$$(2.2) \quad T(z \to z_{k-1}) = \max\{T(z \to z_k - (0,1)), T(z \to z_k - (1,0))\}.$$

We define a random region in the plane as follows. Let $0 > \alpha_1 \geq \alpha_2 \geq \cdots$ and $0 > \beta_1 \geq \beta_2 \geq \cdots$ be nonincreasing integer sequences. Define the set $\Gamma_0 \subset \mathbb{Z}^2$ to be

$$(2.3) \qquad \begin{aligned} &\{(m,k) : m \leq 0, k \leq 0\} \cup \{(m,k) : k > 0, m \leq \alpha_k\} \\ &\cup \{(m,k) : m > 0, k \leq \beta_m\}. \end{aligned}$$

The union of the point $(0,0)$ and the points in the sequences $(\alpha_k, k)_{k>0}$ and $(m, \beta_m)_{m>0}$ is denoted $\gamma_0$. The set $\gamma_0$ coincides with the upper-right corners of the line defining the boundary of $\Gamma_0$; abusing notation, we also call this line $\gamma_0$. See Figure 2. We assume that $\gamma_0$ has asymptotic directions:

$$(2.4) \qquad \lim_{k \to \infty} \frac{\alpha_k}{k} = \frac{-(1-\lambda)}{\lambda} \quad \text{and} \quad \lim_{m \to \infty} \frac{\beta_m}{m} = \frac{-\rho}{1-\rho},$$

for some $\lambda \in (0, 1]$ and $\rho \in [0, 1)$. Then $\gamma_0$ is the boundary of a "cone" containing the complement of $\Gamma_0$, including the positive quadrant, whose asymptotic angle is in $[90°, 270°)$. The angle is in $[90°, 180°)$ if and only if $\lambda > \rho$.

We write $A_k = (\alpha_k + 1, k)$, $k \geq 1$, and $B_m = (m, \beta_m + 1)$, $m \geq 1$. For $z \notin \Gamma_0$, define

$$G^1(z) = \max_{0 < k \leq z(2)} T(A_k \to z)$$

and

$$G^2(z) = \max_{0 < m \leq z(1)} T(B_m \to z),$$

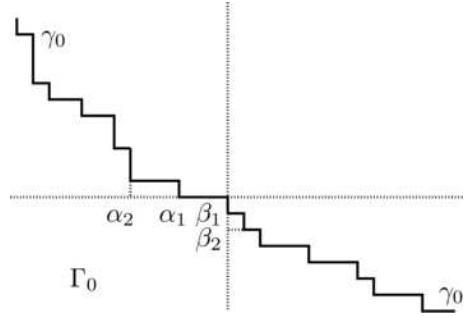

FIG. 2. *The initial configuration: $\alpha_1 = -3$, $\alpha_2 = \alpha_3 = -6$, $\beta_1 = -1$, $\beta_2 = -2$. $\Gamma_0$ is the region below the line $\gamma_0$.*

with, say, $G^1(z) = 0$ if $z(2) \leq 0$ and $G^2(z) = 0$ if $z(1) \leq 0$.

Define

(2.5) $$G(z) = \max\{G^1(z), G^2(z)\}.$$

This is the *last-passage time* from the set

(2.6) $$\{A_k, k \geq 1\} \cup \{B_m, m \geq 1\}$$

to $z$. The quantities $G(z)$ satisfy the recurrence relation

(2.7) $$G(z) = X(z) + \max\{G(z - (0,1)), G(z - (1,0))\},$$

with boundary condition $G(z) = 0$ if $z \in \Gamma_0$. Since the weights are independent with a continuous distribution, with probability 1 all paths have different weights, and so $G(z)$ is achieved by a unique *optimal path* $z_0 \in \Gamma_0$, $z_1, \ldots, z_n = z$, which can be recovered backward with the recurrence relation

(2.8) $$z_{k-1} = \begin{cases} z_k - (1,0), & \text{if } G(z_k - (1,0)) > G(z_k - (0,1)), \\ z_k - (0,1), & \text{otherwise,} \end{cases}$$

so that $z_{k-1} = \arg\max\{G(z_k - (0,1)), G(z_k - (1,0))\}$. Of course, the length $n$ of the optimal path depends on $z$, $\Gamma_0$ and $\mathcal{X}$.

For $t > 0$, define the growth process

(2.9) $$\Gamma_t = \{z : G(z) \leq t\},$$

and the competing growth process

(2.10) $$\Gamma_t^1 = \{z \notin \Gamma_0 : G^2(z) < G^1(z) \leq t\},$$

(2.11) $$\Gamma_t^2 = \{z \notin \Gamma_0 : G^1(z) < G^2(z) \leq t\};$$

see Figure 3. Hence if $z \notin \Gamma_0$ and $G(z) \leq t$, then $z \in \Gamma_t^1$ if the optimal path to $z$ comes from some $A_k$, and $z \in \Gamma_t^2$ if the optimal path to $z$ comes from some $B_m$.

The dynamics of the model can be explained as follows. $\Gamma_t$ is the set of "occupied vertices" at time $t$. We denote by $\gamma_t$ its boundary, the *growth interface* at time $t$. This set is increasing in $t$, and at any $t$ it is a decreasing set: if $z \in \Gamma_t$ and $z' \leq z$, then $z' \in \Gamma_t$ also. Each site $z$ becomes occupied at rate 1, once the sites $z - (0,1)$ and $z - (1,0)$ are both occupied. Each added site joins one of two clusters as follows: sites $(m,k)$ with $m \leq 0, k > 0$ (upper-left quadrant) join cluster 1; sites $(m,k)$ with $m > 0, k \leq 0$ (lower-right quadrant) join cluster 2; sites $(m,k)$ with $m > 0, k > 0$ may join either cluster, depending on which cluster contains sites $(m-1,k)$ and $(m,k-1)$; if these two neighbors both belong to the same cluster, then $(m,k)$ joins this cluster also; if they differ [and then in fact we must have that $(m-1,k)$ is in cluster 1 and $(m,k-1)$ is in cluster 2], then $(m,k)$ joins the cluster of the neighbor which became occupied more recently.

$\Gamma_\infty^j$ is the set of sites which eventually join cluster $j$. $\Gamma_\infty^1$ and $\Gamma_\infty^2$ are connected and have the following "increasing set" properties. If $(m,k) \in \Gamma_\infty^1$, then also $(m,k') \in \Gamma_\infty^1$ for all $k' > k$; if $(m,k) \in \Gamma_\infty^2$, then also $(m',k) \in \Gamma_\infty^2$ for all $m' > m$. The sets $\Gamma_\infty^1$ and $\Gamma_\infty^2$ are separated by the *competition interface*. This consists of the path of sites $z$ such that $z + (0,1) \in \Gamma_\infty^1$ while $z + (1,0) \in \Gamma_\infty^2$. This path starts at $\phi_0 := (0,0)$: we write it as $\phi = (\phi_0, \phi_1, \phi_2, \ldots)$. It can be constructed recursively as follows. Given $\phi_n$, let $\phi_{n+1}$ be equal to $\phi_n + (1,0)$ if $\phi_n + (1,1)$ belongs to $\Gamma_\infty^1$, and be equal to $\phi_n + (0,1)$ if $\phi_n + (1,1)$ belongs to $\Gamma_\infty^2$. Equivalently, $\phi_{n+1}$ equals $\phi_n + (1,0)$ if $G(\phi_n + (1,0)) < G(\phi_n + (0,1))$, and equals $\phi_n + (0,1)$ otherwise [because the point $\phi_n + (1,1)$ takes the color of whichever of its neighbors below or to the left is occupied later]. Thus

(2.12) $$\phi_{n+1} = \arg\min\{G(\phi_n + (1,0)), G(\phi_n + (0,1))\}.$$

So the competition interface starts at $(0,0)$, and thereafter takes steps either to the vertex above or to the vertex to the right, choosing whichever is

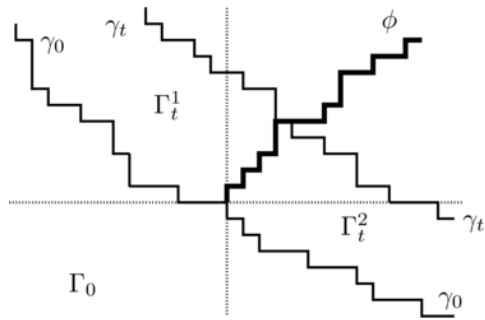

FIG. 3. *Growth and competition interfaces.* $\Gamma_t^1$ *is the set with boundaries* $\gamma_0$, $\gamma_t$ *above* $\phi$ *and* $\Gamma_t^2$ *is the set with boundaries* $\gamma_0$, $\gamma_t$ *below* $\phi$.

occupied first. We show that this interface $\phi$ converges almost surely to an asymptotic direction:

THEOREM 1. *If the initial growth interface $\gamma_0$ satisfies (2.4), then $\mathbb{P}$-a.s. there exists an angle $\theta \in [0, 90°]$ (which may be random, with a distribution depending on $\gamma_0$), such that*

$$\text{(2.13)} \qquad \frac{\phi_n}{|\phi_n|} \to e^{i\theta} := (\cos\theta, \sin\theta) \qquad \text{as } n \to \infty.$$

The convergence of the competition interface is proved in Section 3, using results for point to half-line percolation.

In the case $\lambda > \rho$ especially, we emphasize one particular class of initial distributions of the growth interface. Let $\lambda \in (0, 1]$ and $\rho \in [0, 1)$ and define a random initial growth interface $\gamma_0 = (\gamma_0(j))_{j\in\mathbb{Z}} \subseteq \mathbb{Z}^2$ with $\gamma_0(-1) = (-1, 0)$, $\gamma_0(0) = (0, 0)$, $\gamma_0(1) = (0, -1)$, as follows. Starting from $(-1, 0)$, walk one unit up with probability $\lambda$ and one unit left with probability $1 - \lambda$ to obtain $(\gamma_0(j))_{j<0}$. Then, starting from $(0, -1)$ walk down with probability $\rho$ and right with probability $1 - \rho$ to get $(\gamma_0(j))_{j>0}$. [In this case, the sequences $(\alpha_k - \alpha_{k+1})$ and $(\beta_m - \beta_{m+1})$ are independent, and each is an i.i.d. sequence of geometric random variables taking values in $\{0, 1, 2, \ldots\}$ with means $(1 - \lambda)/\lambda$ and $\rho/(1 - \rho)$, respectively.] We denote by $\nu = \nu_{\lambda,\rho}$ the law of $\gamma_0$.

We show that there is a phase transition from a random to a deterministic direction for the competition interface when the initial growth interface goes from convex to concave:

THEOREM 2. (a) *If $\gamma_0$ satisfies (2.4) with $\lambda \leq \rho$ then the angle $\theta$ in Theorem 1 is $\mathbb{P}$-a.s. constant, and is given by*

$$\text{(2.14)} \qquad \tan\theta = \frac{\lambda\rho}{(1-\lambda)(1-\rho)}.$$

(b) *If $\gamma_0$ satisfies (2.4) with $\lambda > \rho$ then $\mathbb{P}$-a.s the angle $\theta$ satisfies*

$$\text{(2.15)} \qquad \left(\frac{\rho}{1-\rho}\right)^2 \leq \tan\theta \leq \left(\frac{\lambda}{1-\lambda}\right)^2.$$

(c) *If furthermore $\gamma_0$ is distributed according to $\nu_{\lambda,\rho}$ with $\lambda > \rho$, then the distribution of the angle $\theta$ is given by*

$$\text{(2.16)} \qquad \tan\theta = \left(\frac{1-U}{1+U}\right)^2,$$

*where $U$ is a random variable uniformly distributed on $[1 - 2\lambda, 1 - 2\rho]$.*

The condition that $\gamma_0$ is distributed according to the measure $\nu_{\lambda,\rho}$ is essential to the result of (2.16). A local perturbation of the measure would give rise to a different asymptotic law.

In Section 3, we give a direct proof of Theorem 1, which also yields the bounds in (2.15) and the limiting direction in (2.14). The results for $\lambda \leq \rho$ also follow from the law of large numbers in Seppäläinen [26] for the second-class particle and the correspondence of the interface with the second-class particle described later in this paper. This correspondence is also crucial in our proof of (2.16), where we use the weak law of large numbers proved by Ferrari and Kipnis [8] for the path of the second-class particle.

In Section 3 we will also prove and use the following result which describes the linear growth of the passage times in each direction in the positive quadrant:

PROPOSITION 2.1. *Suppose that $\gamma_0$ satisfies* (2.4). *Then with probability* 1,

$$\text{(2.17)} \qquad \lim_{|z| \to \infty, z \in \mathbb{Z}_+^2} \frac{G(z) - p(z)}{|z|} = 0,$$

*where*

$$\text{(2.18)} \quad p(z) = \begin{cases} (\sqrt{z(1)} + \sqrt{z(2)})^2, \\ \qquad \text{if } \left(\frac{\rho}{1-\rho}\right)^2 \leq \frac{z(2)}{z(1)} \leq \left(\frac{\lambda}{1-\lambda}\right)^2, \\ \frac{z(1)}{1-\lambda} + \frac{z(2)}{\lambda}, \\ \qquad \text{if } \frac{z(2)}{z(1)} \geq \max\left(\left(\frac{\lambda}{1-\lambda}\right)^2, \frac{\lambda\rho}{(1-\lambda)(1-\rho)}\right), \\ \frac{z(1)}{1-\rho} + \frac{z(2)}{\rho}, \\ \qquad \text{if } \frac{z(2)}{z(1)} \leq \min\left(\left(\frac{\rho}{1-\rho}\right)^2, \frac{\lambda\rho}{(1-\lambda)(1-\rho)}\right). \end{cases}$$

[The definition in (2.18) makes sense since $\frac{\lambda\rho}{(1-\lambda)(1-\rho)}$ is always between $(\frac{\lambda}{1-\lambda})^2$ and $(\frac{\rho}{1-\rho})^2$, whatever the values of $\lambda$ and $\rho$.]

Note that $p(\alpha z) = \alpha p(z)$ for all $\alpha > 0$ and $z \in \mathbb{Z}_+^2$. Then Proposition 2.1 can easily be rewritten in the form of a "shape theorem," giving an asymptotic shape under linear rescaling for the covered region, or rather its intersection with the positive quadrant. Recall that $\Gamma_t = \{z : G(z) \leq t\}$ and let $\overline{\Gamma}_t = (\Gamma_t + [0,1]^2) \cap \mathbb{R}_+^2$ (we have added a box of area 1 to each point to form a subset of $\mathbb{R}_+^2$ from the set of points in $\mathbb{Z}_+^2$). Now let $\Gamma$ be the set $\{y \in \mathbb{R}_+^2 : p(y) \le 1\}$. Then for all $\varepsilon > 0$, with probability 1, for sufficiently large $t > 0$

(2.19) $$t(1-\varepsilon)\Gamma \subseteq \overline{\Gamma}_t \subseteq t(1+\varepsilon)\Gamma.$$

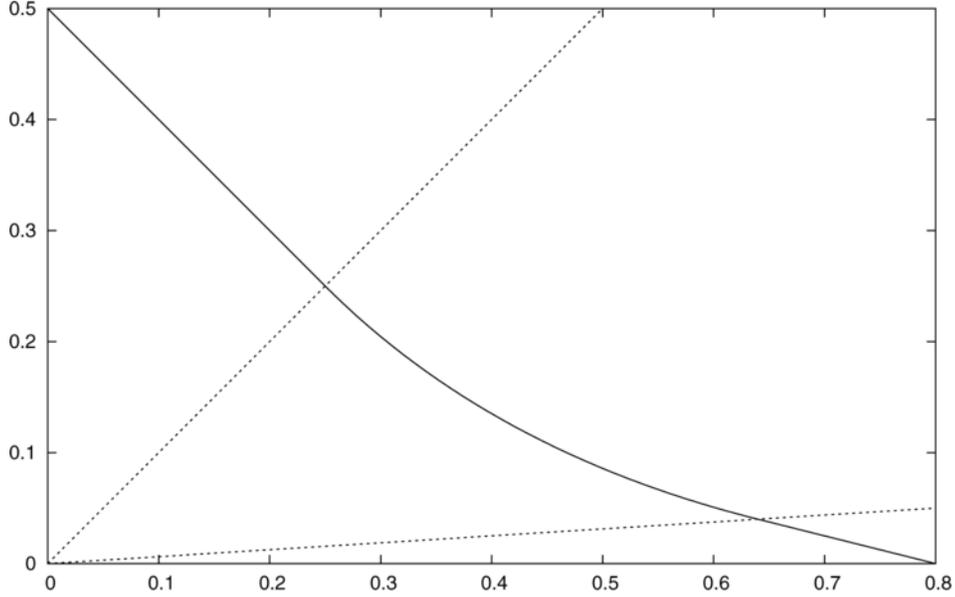

FIG. 4. *A plot of the curve $p(z) = 1$ for $\lambda = 0.5$ and $\rho = 0.2$, together with the lines $\frac{z(2)}{z(1)} = (\frac{\lambda}{1-\lambda})^2$ and $\frac{z(2)}{z(1)} = (\frac{\rho}{1-\rho})^2$ which separate the linear parts of the curve from the curved part when $\rho < \lambda$.*

Equally, this implies that the intersection of the rescaled growth interface $\gamma_t/t$ with the positive quadrant eventually lies between the curves $(1-\varepsilon)\gamma$ and $(1+\varepsilon)\gamma$, where $\gamma$ is the curve $\{y : p(y) = 1\}$. We can relate the law of large numbers for the competition interface to the form of the curve $p(y) = 1$ given in Proposition 2.1. When $\lambda > \rho$, the limiting curve consists of two straight line-segments joined smoothly by a curved part. The angle of the competition interface is random, and its distribution is supported on the cone spanning the curved part. As $\lambda - \rho$ approaches 0 from above, this cone becomes smaller and disappears when $\lambda = \rho$. In this case $p(y) = 1$ is a straight line and the competition interface has a deterministic direction. For $\lambda < \rho$, the curve $p(y) = 1$ consists of two straight line-segments, and the competition interface has a deterministic angle, along the direction of the "shock" joining these two line-segments. See Figure 4.

2.2. *Simple exclusion and second-class particles.* The totally asymmetric simple exclusion process (TASEP) $(\eta_t, t \ge 0)$ is a Markov process in the state space $\{0,1\}^{\mathbb{Z}}$ whose elements are particle configurations. $\eta_t(j) = 1$ indicates a particle at site $j$ at time $t$; otherwise $\eta_t(j) = 0$ (a hole is at site $j$ at time $t$). With rate 1, if there is a particle at site $j$, it attempts to jump to site $j+1$; if there is a hole at $j+1$ the jump occurs, otherwise nothing happens.

To construct a realization of this process à la Harris, one considers independent one dimensional Poisson processes $\mathcal{N} = (N_x(\cdot), x \in \mathbb{Z})$ of intensity 1; let $\mathbb{Q}$ be the law of $\mathcal{N}$. The process $(\eta_t, t \geq 0)$ can be constructed as a deterministic function of the initial configuration $\eta_0$ and the Poisson processes $\mathcal{N}$ as follows: if $s$ is a Poisson epoch of $N_x$ and there is a particle at $x$ and no particle at $x+1$ in the configuration $\eta_{s-}$, then at time $s$ the new configuration is obtained by making the particle jump from $x$ to $x+1$. Let $\Phi$ be the function that takes $\eta_0$ and $\mathcal{N}$ to $(\eta_t, t \geq 0)$.

*Connection to the growth interface.* We relate the simple exclusion process to the growth model by the following method which originates with Rost [25].

We first relate initial configurations for the exclusion processes and initial growth interfaces bijectively: given an initial configuration $\eta_0$, define $\gamma_0$ by $\gamma_0(0) = (0,0)$ and

$$(2.20) \qquad \gamma_0(j) - \gamma_0(j-1) = (1 - \eta_0(j), -\eta_0(j)).$$

Observe that $\gamma_0$ has the asymptotics (2.4) if and only if $\eta_0$ satisfies

$$(2.21) \qquad \lim_{n \to \infty} \frac{1}{n} \sum_{j=-1}^{-n} \eta(j) = \lambda \quad \text{and} \quad \lim_{n \to \infty} \frac{1}{n} \sum_{j=1}^{n} \eta(j) = \rho.$$

We will consider initial configurations $\eta_0$ for which $\eta_0(0) = 0$ and $\eta_0(1) = 1$, so that there is a hole at site 0 and a particle at site 1. This corresponds to the condition that $\gamma_0(0) - \gamma_0(-1) = (1,0)$ and $\gamma_0(1) - \gamma_0(0) = (0,-1)$ (which is what we assumed when setting up the growth model). Furthermore, in the case where $\gamma_0$ is distributed according to the measure $\nu_{\lambda,\rho}$ described after Theorem 1, the distribution of $\eta_0$ is product measure, with density $\lambda$ of particles to the left of site 0 and density $\rho$ of particles to the right of site 1. We will denote this measure by $\nu_{\lambda,\rho}$ also.

Now label the particles sequentially from right to left and the holes from left to right, with the convention that the particle at site 1 and the hole at site 0 are both labeled 0. Let $P_j(0)$ and $H_j(0)$, $j \in \mathbb{Z}$ be the positions of the particles and holes, respectively, at time 0. The position at time $t$ of the $j$th particle $P_j(t)$ and the $i$th hole $H_i(t)$ are functions of the variables $G(z)$ with $z \in C_0 \setminus \gamma_0$ (defined earlier for the growth model) by the following rule: at time $G((i,j))$, the $j$th particle and the $i$th hole interchange positions. Disregarding labels and defining $\eta_t(P_j(t)) = 1$, $\eta_t(H_j(t)) = 0$, $j \in \mathbb{Z}$, the process $\eta_t$ indeed realizes the exclusion dynamics. At time $t$ the particle configuration $\eta_t$ and the growth interface $\gamma_t$ still satisfy the same relation as $\eta_0$ and $\gamma_0$ (2.20).

Note that the shape result in Proposition 2.1 is closely related to the hydrodynamics for the TASEP. The macroscopic density evolution is governed by the Burgers equation; if $\lambda = \rho$, then the density profile is constant, while if the densities to the right and to the left are different, the discontinuity at the origin produces a shock if $\lambda < \rho$ and a rarefaction fan if $\lambda > \rho$. See [10] for further details in this context and for references.

*Coupling and second-class particles.* Let $\eta$ and $\eta'$ be two arbitrary configurations. The *basic* coupling between two exclusion processes with initial configurations $\eta$ and $\eta'$, respectively, is the joint realization $(\Phi(\eta, \mathcal{N}), \Phi(\eta', \mathcal{N})) = ((\eta_t, \eta'_t), t \geq 0)$ obtained by using the same Poisson epochs for the two different initial conditions.

Given a configuration of particles $\eta$, let $\eta'$ be a configuration which differs from $\eta$ only at the origin. Call $X(0) = 0$ the site where both configurations differ at time zero. With the basic coupling, the configurations at time $t$ differ only at the site $X(t)$ defined by

$$(2.22) \qquad X(t) := \sum_x x \mathbf{1}\{\eta_t(x) \neq \eta'_t(x)\}.$$

$(X(t), t \geq 0)$ is the trajectory of a "second-class particle." The process $((\eta_t, X(t)), t \geq 0)$ is Markovian but the process $(X(t), t \geq 0)$ is not. The motion of $X(t)$ depends on the configuration of $\eta'_t$ in its neighboring sites. The second-class particle jumps one unit to the right at rate 1 if there is no $\eta'$ particle in its right nearest neighbor and it jumps one unit to the left at rate 1 if there is a $\eta'$ particle in its left nearest neighbor site, interchanging positions with it.

Define the process $\psi(t)$ by

$$(2.23) \qquad \psi(t) := (I(t), J(t)) := \phi_n \qquad \text{for } t \in [G(\phi_n), G(\phi_{n+1})).$$

Note that $\psi(t)$ gives the position of the farthest (to the North-East) intersecting point between the competition interface $\phi$ and the growth interface $\gamma_t$ (Figure 5). The relation between competition interfaces and second-class particles is given by the following proposition:

PROPOSITION 2.2. *There exists a coupling of the growth process $(\gamma_t, t \geq 0)$ and the exclusion process with second-class particle $((\eta_t, X(t)), t \geq 0)$ under which, for all $t$, $J(t)$ equals the number of leftward jumps of the second-class particle in $[0, t]$, and $I(t)$ is the number of rightward jumps of the second-class particle in $[0, t]$. Thus the trajectory $(X(t), t \geq 0)$ is identical to the trajectory $(I(t) - J(t), t \geq 0)$.*

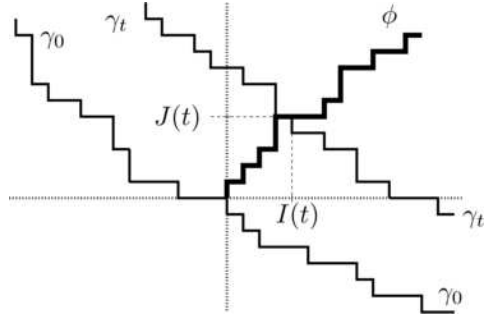

Fig. 5. *The position of the interface at time t.*

OUTLINE OF PROOF. The proof of this proposition can be carried out using the techniques introduced by Ferrari and Pimentel in [11]. There are two main elements. First, the exclusion process with a second-class particle is shown to be equivalent to an exclusion process with no second-class particle, but with an extra site. In this construction, the second-class particle is replaced by a pair of sites, of which the left site contains a hole and the right site contains a particle. When this is done, an initial configuration with a second-class particle at the origin becomes an initial configuration with a hole at the origin and a particle at site 1. The second element is then the correspondence described around (2.20) between the exclusion process (specifically, one starting from such an initial configuration) and the growth process. For details of the construction, see Proposition 3 and Lemma 6 of Ferrari and Pimentel [11]. □

We obtain a strong law of large numbers for the second-class particle:

THEOREM 3. *Suppose the initial condition $\eta_0$ satisfies* (2.21).

(a) *The limit*

$$U = \lim_{t \to \infty} \frac{X(t)}{t} \tag{2.24}$$

*exists $\mathbb{Q}$-a.s.*

(b) *If $\lambda > \rho$, then $U \in [1 - 2\lambda, 1 - 2\rho]$ a.s.*

(c) *If moreover $\eta_0$ is distributed according to the product measure $\nu_{\lambda,\rho}$ with $\lambda > \rho$, then $U$ has the uniform distribution on $[1 - 2\lambda, 1 - 2\rho]$.*

(d) *If $\lambda \leq \rho$, then $U = 1 - \lambda - \rho$ a.s.*

We give an independent proof of the strong law of the almost sure convergence for all $\lambda$ and $\rho$, using the correspondence between the second-class particle and the competition interface. This approach also gives the bounds on the distribution in part (b) [and could also easily be used to provide the limiting value for $\lambda < \rho$ in part (d)]. We include the statements in parts (c) and (d) for completeness. The convergence and the limit in part (d) for the case $\lambda \leq \rho$ was proved by Ferrari [5] in the case of product measure $\nu_{\lambda,\rho}$ and by Seppäläinen [26] in the general case of initial configurations satisfying (2.21). Rezakhanlou [24] proved, convergence in probability for more general initial conditions when the limit is not random. In the case of part (c) of product measure $\nu_{\lambda,\rho}$ with $\lambda > \rho$, Ferrari and Kipnis [8] proved convergence in distribution to the uniform distribution, and Mountford and Guiol [21] proved the almost sure convergence.

In Section 3, we give the proofs of Theorem 1 and of parts (a) and (b) of Theorem 2, along with Proposition 2.1. Using these results, we can prove the remaining properties stated above:

PROOFS OF THEOREMS 3(a), (b) AND 2(c). From Theorem 1, we know that with probability 1, the limit $I(t)/J(t) \to \tan \theta$ exists for some (maybe random) angle $\theta$.

Also $(I(t), J(t)) = \psi(t) \in \gamma_t$, and so by the shape theorem (Proposition 2.1), we have that $p(I(t)/t, J(t)/t) = t^{-1}p(\psi(t)) \to 1$ with probability 1, as $t \to \infty$.

Hence $(I(t)/t, J(t)/t) \to (I, J)$ as $t \to \infty$, where $(I, J)$ is the point on the curve $p(y) = 1$ which is at angle $\theta$ from the origin. Since by Proposition 2.2 we have $X(t) = I(t) - J(t)$, this implies that $X(t)/t$ also converges with probability 1, as required for Theorem 3(a). Its limit, $U$ say, is given by $U = I - J$.

If $\lambda > \rho$, we have from Theorem 2(b) that the limiting direction $\theta$ satisfies (2.15). In this case $p(I, J) = 1$ implies that $\sqrt{I} + \sqrt{J} = 1$. Putting this together with $U = I - J$, we obtain $I = (1 + U)^2/4$ and $J = (1 - U)^2/4$. Substituting $\tan \theta = J/I$ and using (2.15), we obtain $U \in [1 - 2\lambda, 1 - 2\rho]$ as required for part (b) of Theorem 3.

To prove Theorem 2(c), we will use the weak law of large numbers for the second-class particle in the case of an initial configuration $\nu_{\lambda,\rho}$ with $\lambda > \rho$ (from Ferrari and Kipnis [8]), which gives that $X(t)/t$ converges in distribution to the uniform distribution on $[1 - 2\lambda, 1 - 2\rho]$. So this is also the distribution of the almost sure limit $U$. Then from the previous paragraph we obtain

$$\tan \theta = \frac{J}{I} = \left(\frac{1 - U}{1 + U}\right)^2, \tag{2.25}$$

as required. □

We have shown a law of large numbers for competition interfaces and second-class particles under conditions (2.4) and (2.21), respectively. However, for $\lambda > \rho$, we can compute the law of the asymptotic direction only when $\gamma$ and $\eta$ have distribution $\nu_{\lambda,\rho}$. It would be nice to have a description of those laws in more general cases. In this direction, we can say that it depends strongly on the microscopic details of the initial condition, and we have an example which shows such dependence. Let $L > 0$ and consider an initial growth interface $\gamma^L$ where $\alpha_k = -L$ and $\beta_k = L$ for all $k \geq 1$. In the exclusion context, this corresponds to an initial configuration $\eta^L$ with $\eta^L(k) = 0$ for $k = -1, \ldots, -L$ and $k > L$ and with $\eta^L(k) = 1$ for $k = 0, \ldots, L$ and $k < -L$. Thus for all $L > 0$, $\gamma^L$ satisfies (2.4) as well as $\eta^L$ satisfies (2.21), with $\lambda = 1$ and $\rho = 0$. By Theorems 1 and 3, there exist $\theta^L$ and $U^L$ which correspond to the asymptotic direction of the competition interface and of the second-class particle, respectively. Using the same approach as in [23], one can show that $\theta^L \to \pi/4$ in probability as $L \to \infty$ (this would correspond to a multi-type shape theorem for the LPP model). Correspondingly for the second class particle one has that $U^L \to 0$ in probability. Thus we see that the law of $\theta^L$ is quite different for $L = 1$ and for large $L$. More generally it is not hard to see that the laws of $\theta$ or $U$ are sensitive to local changes in the initial condition. For example, consider starting from the distribution $\nu_{1,0}$ in which with probability 1 every negative site has a particle and every positive site has a hole. In this case, the state at any finite time $t$ differs only at finitely many sites from the initial state; nevertheless $\theta$ and $U$ converge almost surely to nondeterministic limits. Such convergence is only possible if their distributions are sensitive to local changes.

### 2.3. Fluctuations.

*Concave sector: central limit theorem.* We now turn to questions concerning the fluctuation of the competition interface around its asymptotic direction (random or deterministic).

When $\lambda < \rho$, the cone into which the occupied set grows has a concave angle. In the case of an initial configuration distributed according to the product measure $\nu_{\lambda,\rho}$, the interface position on the scale $\sqrt{t}$ is determined by the initial configuration. This makes it possible to compute explicitly the joint asymptotic law of $(I(t), J(t))$ after centering and rescaling by $\sqrt{t}$. The proof is based on results in the same vein for the second-class particle and the flux of particles. Let

$$(2.26) \qquad N^\eta_{r,r'} = \sum_{i=r}^{r'} \eta(i)$$

with the convention $\sum_r^{r'} = -\sum_{r'}^{r}$ if $r' < r$. Denote $I^\eta(t)$ and $J^\eta(t)$ the coordinates of $\psi(t)$ when the initial particle configuration is $\eta$. Denote $1 - \eta$ the configuration defined by $(1 - \eta)(x) = 1 - \eta(x)$.

THEOREM 4. *Dependence of the initial configuration and the central limit theorem:*

$$\lim_{t\to\infty}\frac{1}{t}\int\nu_{\lambda,\rho}(d\eta)\mathbb{E}\bigg(I^\eta(t)-\frac{(1-\lambda)N^{1-\eta}_{0,(\rho-\lambda)t}-(1-\rho)N^\eta_{-(\rho-\lambda)t,0}}{\rho-\lambda}$$

$$-\lambda(1-\rho)t\bigg)^2 \quad (2.27)$$

$$=0,$$

$$\lim_{t\to\infty}\frac{1}{t}\int\nu_{\lambda,\rho}(d\eta)\mathbb{E}\bigg(J^\eta(t)-\frac{-\lambda N^{1-\eta}_{0,(\rho-\lambda)t}+\rho N^\eta_{-(\rho-\lambda)t,0}}{\rho-\lambda}-\lambda(1-\rho)t\bigg)^2$$

$$(2.28) \quad =0.$$

*The limiting covariance matrix is given by*

$$(2.29) \qquad \lim_{t\to\infty}\frac{1}{t}\mathbb{V}I(t)=\frac{(1-\rho)(1-\lambda)(\rho(1-\lambda)+\lambda(1-\rho))}{\rho-\lambda},$$

$$(2.30) \qquad \lim_{t\to\infty}\frac{1}{t}\mathbb{V}J(t)=\frac{\lambda\rho(\rho(1-\lambda)+\lambda(1-\rho))}{\rho-\lambda},$$

$$(2.31) \qquad \lim_{t\to\infty}\frac{1}{t}\mathrm{Cov}(I(t),J(t))=-\frac{2\lambda(1-\lambda)\rho(1-\rho)}{\rho-\lambda},$$

*and the central limit theorem holds: the vector*

$$(2.32) \qquad \frac{1}{\sqrt{t}}(I(t)-(1-\rho)(1-\lambda)t,J(t)-\rho\lambda t)$$

*converges in law to a two-dimensional normal distribution with covariance matrix* (2.29), (2.30) *and* (2.31), *as* $t\to\infty$.

*Fluctuations: convex or flat sectors.* In the case $\lambda\geq\rho$, we relate the competition interface with the semi-infinite geodesics defined at the beginning of the section.

Ferrari and Pimentel [11] proved that the following holds with probability 1: for every semi-infinite geodesic $(z_k)_{k>0}$, there exists $\alpha\in[0,\pi/2]$ such that $(z_k)_{k>0}$ has asymptotic direction $e^{i\alpha}$, and for every $\alpha\in[0,\pi/2]$ there exists a semi-infinite geodesic with asymptotic direction $e^{i\alpha}$.

For $r>0$, $\alpha\in(0,\pi/2)$ and $z\in\mathbb{Z}^2$, let $\mathbb{H}_z(\alpha,r)$ be the hyperplane that is perpendicular to the angle $\alpha$ and contains $z+re^{i\alpha}$. Given a semi-infinite up-right path $(z_k)_{k>0}$ with asymptotic direction $e^{i\alpha}$, we denote by $y_{z_1}(\alpha,r)$ the intersection between $\mathbb{H}_{z_1}(\alpha,r)$ and $(z_k)_{k>0}$.

We say that a path $(z_k)_{k>0}$, with direction $e^{i\alpha}$, "has fluctuation at most $\kappa$" if $|y_{z_0}(\alpha,r)-(z_1+re^{i\alpha})|\leq r^\kappa$ for all sufficiently large $r$.

Now define fluctuation exponents of the set of geodesics by

$$\xi_1(\alpha) = \inf\{\kappa : \mathbb{P}(\text{all semi-infinite geodesics with}$$

$$\text{direction } e^{i\alpha} \text{ have fluctuation at most } \kappa) = 1\},$$

for $\alpha \in (0, \pi/2)$, and

(2.33)  $\xi_2 = \inf\{\kappa : \mathbb{P}(\text{all semi-infinite geodesics have fluctuation at most} \kappa) = 1\}.$

In fact, Balázs, Cator and Seppäläinen [1] prove probability estimates giving upper and lower bounds on the fluctuations of finite geodesics which can be used to show that $\xi_1(\alpha) = 2/3$ for all $\alpha \in (0, \pi/2)$.

It is immediate that $\xi_2 \geq \sup_\alpha \xi_1(\alpha)$, but the reverse inequality is not obvious (since we consider uncountably many directions $\alpha$). The method used to prove the existence of semi-infinite geodesics (see [11]) yields the upper bound $\xi_2 \leq 3/4$. However, under an assumption on the asymptotic behavior of $\mathbb{P}(T(z' \to z) \leq n)$, as in (1.5) of Johansson [16] (see also the last paragraph of Section 1 in the same work), one can mimic the argument developed by Wüthrich [28] to prove that

$$\xi_1(\alpha) \leq \xi_2 \leq 2/3, \tag{2.34}$$

which, together with the bounds from [1], would imply $\xi_2 = 2/3$.

In an analogous way we define the fluctuation exponent $\chi$ of the competition interface by

(2.35)  $\chi = \inf\{\kappa : \mathbb{P}(\text{the competition interface has fluctuation at most } \kappa) = 1\}.$

(Of course the exponent $\chi$ depends on the initial configuration $\gamma_0$.)

In equilibrium we exhibit a duality between competition interfaces and semi-infinite geodesics and prove that they have the same fluctuation exponent (which, from the discussion above, is equal to 2/3). This duality is also developed in [1].

THEOREM 5. *Assume that $\gamma_0$ has distribution $\nu_{\lambda,\rho}$ with $\lambda = \rho$. Then $\chi = \xi_1(\alpha)$, where $\tan(\alpha) = \rho^2/(1-\rho)^2$.*

For a convex initial growth interface ($\lambda > \rho$) we show that the fluctuations of semi-infinite geodesics dominate the fluctuations of competition interfaces:

THEOREM 6. *Assume that $\lambda > \rho$ and that $\gamma_0$ satisfies* (2.4). *Then $\chi \leq \xi_2$.*

This is proved by showing that the competition interface (with random direction $\theta$) is contained between two disjoint semi-infinite geodesics with direction $\theta$. Observe that this implies that the random direction $\theta$ chosen by the interface is somehow exceptional, since it can be shown that for any fixed direction $\alpha$, the probability that there exist two disjoint geodesics with direction $\alpha$ is 0 [although with probability 1, the set of such exceptional directions is dense in $(0, \pi/2)$]. It would be natural to conjecture that in fact $\chi = 2/3$ in this case also.

The central limit theorem in the last section for $\lambda < \rho$ with initial distribution $\nu_{\lambda,\rho}$ indicates that in that case $\chi = 1/2$. However, as already observed, the fluctuations in that case are controlled by the fluctuations of the initial growth interface. If, for example, the growth interface was replaced by a "periodic" one with the same asymptotic directions, one would expect to see fluctuations of order $1/3$ (as observed by Derrida and Dickmann [4] in an analogous first-passage percolation model).

Theorems 4–6 are proved in Section 4.

Because of the correspondence between the competition interface and the second-class particle (Proposition 2.2) one would expect that for all initial conditions, the competition interface and the second-class particle have the same order of fluctuations. In order to derive this in general, one would need to control the order of the fluctuations of the boundary of the occupied region, for which the first-order behavior is given by the shape theorem (Proposition 2.1).

**3. Point to half-line geodesics.** Consider the point to half-line percolation model as follows: define $\{B_m, m \geq 1\}$, where $B_m = (m, \beta_m + 1)$, and $\beta_1, \beta_2, \ldots$ is a nonincreasing sequence taking negative integer values. Assume that

$$\lim_{m \to \infty} \frac{-\beta_m}{m} = \frac{\rho}{1 - \rho} := d_\rho. \tag{3.1}$$

We have

$$G^2(z) = \max_{0 < m \leq z(1)} \{T(B_m \to z)\}. \tag{3.2}$$

In this section we will develop the properties of the point to half-line percolation model in order to prove Theorem 1 and parts (a) and (b) of Theorem 2, along with the shape result of Proposition 2.1.

Let $M(z)$ be the $m$ which maximizes in the expression above (i.e., the horizontal coordinate of the point from which the maximizing path to $z$ starts):

$$M(z) = \operatorname*{arg\,max}_{0 < m \leq z(1)} \{T(B_m \to z)\}. \tag{3.3}$$

We note the importance of the "critical direction" $d_\rho^2$. Its significance is as follows. Suppose $|z_i| \to \infty$ with $z_i(2)/z_i(1) \to w$. We will show that if $w > d_\rho^2$, the contact point $M(z_i)$ stays small [in fact, one can show by extending the methods below that, with probability 1, $M(z_i)$ is eventually constant]. On the other hand, if $w < d_\rho^2$, then $M(z_i)$ grows linearly with $|z_i|$.

The first aim is to establish the following proposition, which is the key to showing Theorem 1 in the case $\lambda > \rho$.

PROPOSITION 3.1. *Let $\tilde{w} > w > d_\rho^2$. Let $|x_i| \to \infty$ with $x_i(2)/x_i(1) \to w$. Let $|y_j| \to \infty$ with $y_j(2)/y_j(1) \to \tilde{w}$. Then there are a.s. only finitely many $i$ such that, for some $j$, the optimal path from $\{B_m, m \geq 1\}$ to $y_j$ passes through $x_i$.*

3.1. *Sublinear growth.* For $z \geq 0$, let $\mu(z) = (\sqrt{z(1)} + \sqrt{z(2)})^2$. The function $\mu$ identified by Rost [25] represents the "asymptotic shape" for the last-passage percolation model. One has for example that

$$\lim_{|z| \to \infty} \frac{\mathbb{E}T((0,0) \to z)}{\mu(z)} = 1. \tag{3.4}$$

The following moderate deviations estimate is useful several times (proved in Lemma 12 of [11]):

LEMMA 3.1. *For all $\eta > 0$ there is a $b > 0$ such that for all $z_1 \leq z_2$ and for all $r \in [|z_2 - z_1|^{1/2+\eta}, |z_2 - z_1|^{3/2-\eta}]$,*

$$(3.5) \quad \mathbb{P}(|T(z_1 \to z_2) - \mu(z_2 - z_1)| \geq r) \leq b \exp(-br|z_2 - z_1|^{-1/2}).$$

To prove Proposition 3.1, first we show that the growth of $M(z_i)$ is sublinear:

PROPOSITION 3.2. *If $|z_i| \to \infty$ with $z_i(2)/z_i(1) \to w > d_\rho^2$, then $M(z_i)/|z_i| \to 0$ a.s.*

To prove this we need a couple of estimates:

LEMMA 3.2. *Let $\delta, \varepsilon > 0$. Then there exist $\eta > 0$, $c > 0$ and $R < \infty$ such that if $z(2)/z(1) > d_\rho^2 + \delta$ with $n > \varepsilon|z|$ and $|z| > R$, then*

$$\mu(z - (0,1)) - \mu(z - (n, -[d_\rho + \eta]n)) > c|z|. \tag{3.6}$$

PROOF. Fix $z$ and let $\eta > 0$. For $0 \leq x < z(1)$, let
$$f(x) = \mu(z - (x, -[d_\rho + \eta]x))$$
$$= (\sqrt{z(1) - x} + \sqrt{z(2) + (d_\rho + \eta)x})^2$$
$$= g(x)^2,$$

say. We have $f'(x) = 2g(x)g'(x)$.

Suppose that $\eta$ is small enough that

(3.7) $$\frac{d_\rho + \eta}{(d_\rho^2 + \delta)^{1/2}} \leq 1 - u,$$

for some $u > 0$. Then, using $z(2) > (d_\rho^2 + \delta)z(1)$, we have that for all $x \geq 0$,

$$\begin{aligned}
g'(x) &= -\tfrac{1}{2}\{[z(1) - x]^{-1/2} - [d_\rho + \eta][z(2) + (d_\rho + \eta)x]^{-1/2}\} \\
&\leq -\tfrac{1}{2}\{z(1)^{-1/2} - [d_\rho + \eta]z(2)^{-1/2}\} \\
&\leq -\tfrac{1}{2}\{z(1)^{-1/2} - [1 - u]z(1)^{-1/2}\} \\
&= -uz(1)^{-1/2}/2.
\end{aligned}$$

Also for all $x \geq 0$, $g(x) \geq \sqrt{z(2)} \geq d_\rho\sqrt{z(1)}$. So $f'(x) = 2g(x)g'(x) \leq -ud_\rho$ for all $x \geq 0$. Then

$$\mu(z) - \mu(z - (n, -[d_\rho + \eta]n)) = f(0) - f(n)$$

(3.8) $$\geq ud_\rho n$$

$$\geq \varepsilon u d_\rho |z| \qquad \text{for } n \geq \varepsilon|z|.$$

Finally we need to bound $\mu(z) - \mu(z - (0, 1))$:

$$\begin{aligned}
\mu(z) - \mu(z - (0, 1)) &= (\sqrt{z(1)} + \sqrt{z(2)})^2 - (\sqrt{z(1)} + \sqrt{z(2) - 1})^2 \\
&= 1 + 2\sqrt{z(1)}(\sqrt{z(2)} - \sqrt{z(2) - 1}),
\end{aligned}$$

which is less than some constant $\tilde{c}$, uniformly over $z$ satisfying $z(2) > d_\rho^2 z(1)$.

Combining this with (3.8),

(3.9) $$\mu(z - (0, 1)) - \mu(z - (n, -[d_\rho + \eta]n)) > \varepsilon u d_\rho |z| - \tilde{c} > c|z|$$

for all large enough $|z|$ and suitable constant $c$, as required. $\square$

LEMMA 3.3. *Let $\varepsilon, \delta > 0$. Then there exists $\tilde{b}$ such that if $|z|$ is sufficiently large and $z(2)/z(1) > d_\rho^2 + \delta$, then*

(3.10) $$\mathbb{P}(M(z) > \varepsilon|z|) \leq |z|\tilde{b}\exp(-\tilde{b}|z|^{1/2}).$$

PROOF.

(3.11) $$\begin{aligned}
&\mathbb{P}(M(z) > \varepsilon|z|) \\
&\leq \sum_{\varepsilon|z| < m \leq z(1)} \mathbb{P}\{T((m, \beta_m + 1) \to z) > T((1, \beta_1 + 1) \to z)\}.
\end{aligned}$$

Choose $R$ and $\eta$ according to Lemma 3.2. Assume that $|z| > R$, and that $|z|$ is large enough that for all $m \geq \varepsilon|z|$, we have $\beta_m > -(d_\rho + \eta)m$. Note also that $\beta_1 + 1 \leq 0$. Using Lemma 3.2, we bound each term in the sum above:

$$\mathbb{P}\{T((m, \beta_m + 1) \to z) > T((1, \beta_1 + 1) \to z)\}$$
$$\leq \mathbb{P}\{T((m, -[d_\rho + \eta]m) \to z) > T((1, 0) \to z)\}$$
$$= \mathbb{P}\{[T((m, -[d_\rho + \eta]m) \to z) - \mu(z - (m, -[d_\rho + \eta]m))]$$
$$\quad - [T((1, 0) \to z) - \mu(z - (1, 0))]$$
$$\quad > \mu(z - (1, 0)) - \mu(z - (m, -[d_\rho + \eta]m))\}$$
$$\leq \mathbb{P}\{[T((m, -[d_\rho + \eta]m) \to z) - \mu(z - (m, -[d_\rho + \eta]m))]$$
$$\quad - [T((1, 0) \to z) - \mu(z - (1, 0))] > c|z|\}$$
$$\leq \mathbb{P}\{|T(z_1) - \mu(z_1)| > c|z|/2\} + \mathbb{P}\{|T(z_2) - \mu(z_2)| > c|z|/2\},$$

where $z_1 = z - (1, 0)$ and $z_2 = z - (m, -[d_\rho + \eta]m)$. We have $|z|/2 < |z_1| < |z|$, and $v|z| \leq z(2) \leq |z_2| \leq (1 + d_\rho + \eta)|z|$, where $v = (d_\rho^2 + \delta)/(d_\rho^2 + \delta + 1)$. Thus for large enough $|z|$, we can apply the moderate deviations estimate from Lemma 3.1 to bound the last expression above by $\tilde{b}\exp(-\tilde{b}|z|^{1/2})$. Summing over $m$ then gives the desired result. □

PROOF OF PROPOSITION 3.2. Using the bound in Lemma 3.3 and summing over $i$, we have that $\sum_i \mathbb{P}(M(z_i) > \varepsilon|z_i|) < \infty$ for any $\varepsilon$. Applying Borel–Cantelli then gives the result. □

3.2. *Good pairs and good points.* We need some definitions and results which extend the idea of "$h$-straightness" (from Newman [22] and later papers).

For $x \leq z \in \mathbb{Z}^2$, let $R^{\text{out}}[x, z]$ be the set of $z' \geq z$ such that the optimal path from $x$ to $z'$ passes through the point $z$.

For $x \leq z \in \mathbb{Z}^2$ and $\phi > 0$, define the cone $\mathcal{C}(x, z, \phi)$ to be the set of those $z' \geq x$ such that $\text{angle}(z - x, z' - x) \leq \phi$.

Fix $\varepsilon > 0$ and $0 < \delta < 1/4$. For $x \leq z$, say that $(x, z)$ is a *good pair* if

$$(3.12) \quad R^{\text{out}}[x, z] \subseteq \mathcal{C}(x, z, |z - x|^{-\delta}) \cup \bigcup_{\substack{z' \in \mathcal{C}(x, z, |z-x|^{-\delta}) \\ |z' - x| > 2|z - x|}} R^{\text{out}}[x, z'].$$

(This is almost the same thing as saying that

$$(3.13) \quad R^{\text{out}}[x, z] \subseteq \mathcal{C}(x, z, |z - x|^{-\delta}) \cup \{z' : |z' - x| > 2|z - x|\},$$

except for some possible differences at the boundary of the cone.)

Say that $z > 0$ is a *good point* if $(B_m, z)$ is a good pair for every $0 \leq m \leq \varepsilon|z|$.

Let $d(x,y)$ be the Euclidean distance between two points $x$ and $y$, and $d(x,A)$ the minimal Euclidean distance between a point $x$ and a set $A$. Let $[x,y]$ be the straight line segment joining $x$ and $y$.

We use the following lemma, which is a direct consequence of Lemma 11 of [11]. (The proof of that lemma uses the moderate deviations estimate in Lemma 3.1 and an estimate on the convexity of $\mu$ which is essentially Lemma 2.1 of Wüthrich [28].)

LEMMA 3.4. *Let $0 < \delta' < 1/4$. Then there exist positive and finite $C_1$, $C_2$, $C_3$, $C_4$ such that if $|z'| \geq C_1$ and $d(z, [(0,0), z']) \geq |z'|^{1-\delta'}$, then*

$$(3.14) \qquad \mathbb{P}(z' \in R^{\mathrm{out}}[(0,0), z]) \leq C_2 \exp(-C_3 |z'|^{C_4}).$$

LEMMA 3.5. *There exist positive and finite $\tilde{C}_1$, $\tilde{C}_2$ such that for $|z| \geq \tilde{C}_1$,*

$$(3.15) \qquad \mathbb{P}(((0,0), z) \text{ is not a good pair}) \leq \tilde{C}_2 |z|^2 \exp(-C_3 |z|^{C_4}).$$

PROOF. From the definition of a good pair, the following is sufficient for $((0,0), z)$ to be a good pair: for all $z'$ with $|z| \leq |z'| \leq 3|z|$ and $z' \notin \mathcal{C}((0,0), z, |z|^{-\delta})$, we have $z' \notin R^{\mathrm{out}}[(0,0), z]$.

We bound the probability that this fails for any such $z'$. So, suppose that $|z| \leq |z'| \leq 3|z|$ and $\mathrm{angle}(z, z') \geq |z|^{-\delta}$. Choose $\delta'$ with $\delta < \delta' < 1/4$ and then let $\tilde{C}_1$ be sufficiently large that:

(i) $\tilde{C}_1 \geq C_1$;
(ii) if $r \geq \tilde{C}_1$, then $\frac{1}{2}(\frac{r}{3})^{1-\delta} \geq r^{1-\delta'}$;
(iii) if $r \geq \tilde{C}_1$, then $\sin(r^{-\delta}) \geq r^{-\delta}/2$.

Assume that $|z| > \tilde{C}_1$. Then

$$d(z, [(0,0), z']) \geq |z| \sin(\mathrm{angle}(z, z')) \geq |z| \sin(|z|^{-\delta})$$

$$\geq \frac{1}{2}|z|^{1-\delta} \geq \frac{1}{2}\left|\frac{z'}{3}\right|^{1-\delta} \geq |z'|^{1-\delta'}.$$

Then from Lemma 3.4,

$$\mathbb{P}(z' \in R^{\mathrm{out}}[(0,0), z]) \leq C_2 \exp(-C_3 |z'|^{C_4})$$
$$\leq C_2 \exp(-C_3 |z|^{C_4})$$

for all such $z'$. Since there are fewer than $9|z|^2$ such $z'$, we can sum to get the result. □

PROPOSITION 3.3. *Let $\varepsilon > 0$, $0 < \delta < 1/4$. With probability 1, all but finitely many $z > 0$ are good points.*

PROOF. As before, choose $K$ large enough that $|B_m| \leq Km$ for all $m \geq 1$. Suppose that $|z| \geq (1-\varepsilon)^{-1}\tilde{C}_1$. If $m \leq \varepsilon|z|$, then $\tilde{C}_1 \leq (1-\varepsilon)|z| \leq |z - B_m| \leq (K+1)|z|$. By translation invariance, one has that

(3.16) $\quad \mathbb{P}((x,z) \text{ is a good pair}) = \mathbb{P}(((0,0), z-x) \text{ is a good pair}).$

So, using Lemma 3.5, we have

$$\mathbb{P}(z \text{ is not a good point}) \leq \sum_{1 \leq m \leq \varepsilon|z|} \mathbb{P}((B_m, z) \text{ is not a good pair})$$

$$= \sum_{1 \leq m \leq \varepsilon|z|} \mathbb{P}(((0,0), z - B_m) \text{ is not a good pair})$$

$$\leq \varepsilon|z|\tilde{C}_2[(K+1)|z|]^2 \exp(-C_3[(1-\varepsilon)|z|]^{C_4}).$$

This sums to a finite amount over all $z \in \mathbb{Z}_+^2$, so the proposition follows from Borel–Cantelli. $\square$

COROLLARY 3.1. *For any $0 < \varepsilon < 1/2$ and $0 < \delta < 1/4$, there is $C = C(\delta)$ such that, with probability 1, for all but finitely many $z > 0$ one has*

(3.17) $\quad\quad\quad\quad\quad R^{\text{out}}[B_m, z] \subseteq \mathcal{C}(B_m, z, C|z|^{-\delta})$

*for all $m \leq \varepsilon|z|$.*

PROOF. The idea is the same as for the end of the proof of Proposition 3.2 in [22].

From Proposition 3.3, there is a.s. some $L$ such that every $z > 0$ with $|z| \geq L$ is a good point.

Choose any such $z$ and let $x = B_m$ for some $1 \leq m \leq \varepsilon|z|$. The "good point" property gives that

(3.18) $\quad R^{\text{out}}[x, z] \subseteq \mathcal{C}(x, z, |z-x|^{-\delta}) \cup \bigcup_{\substack{z' \in \mathcal{C}(x,z,|z-x|^{-\delta}) \\ |z'-x|>2|z-x|}} R^{\text{out}}[x, z'].$

Note that if $z' \in \mathcal{C}(x, z, \phi_1)$ and $z'' \in \mathcal{C}(x, z', \phi_2)$, then $z'' \in \mathcal{C}(x, z, \phi_1 + \phi_2)$. So, applying the good point property repeatedly, we can obtain by induction that for $m = 1, 2, \ldots,$

(3.19) $\quad\quad R^{\text{out}}[x, z] \subseteq \mathcal{C}(x, z, \eta^{(m)}) \cup \bigcup_{\substack{\tilde{z} \in \mathcal{C}(x,z,\eta^{(m)}) \\ |\tilde{z}-x|>2^m|z-x|}} R^{\text{out}}[x, \tilde{z}],$

where

$$\eta^{(m)} = \sum_{k=0}^{m-1} (2^k|z-x|)^{-\delta} \leq (1-2^{-\delta})^{-1}|z-x|^{-\delta} \leq C(\delta)|z|^{-\delta},$$

say, since $|z - x| \geq (1 - \varepsilon)|z| \geq |z|/2$.

Now the intersection of $R^{\text{out}}[x, \tilde{z}]$ with any finite subset of $\mathbb{Z}_+^2$ is eventually empty as $|\tilde{z} - x| \to \infty$, so we can let $m \to \infty$ to obtain that $R^{\text{out}}[x, z] \subseteq \mathcal{C}(x, z, C|z|^{-\delta})$ as desired. $\square$

PROOF OF PROPOSITION 3.1. Put another way: there are only finitely many points $x_i$ in the sequence such that $y_j \in R^{\text{out}}[B_{M(x_i)}, x_i]$ for some $j$.

We know from Proposition 3.2 that $M(x_i)/|x_i| \to 0$ a.s. So it is enough to show that for some $\varepsilon > 0$, there are a.s. only finitely many points $i$ such that $y_j \in R^{\text{out}}[B_m, x_i]$ for some $j$ and some $0 \leq m \leq \varepsilon|x_i|$.

Let $\eta$ be small enough that $\tilde{w} - \eta > w + \eta$. For large enough $j$ we have $y_j(2)/y_j(1) > \tilde{w} - \eta$, and for large enough $i$ we have $x_i(2)/x_i(1) < w + \eta$.

There exists $K < \infty$ such that, for all $m \geq 1$, $|B_m| \leq Km$. Hence one can find $\varepsilon$ and $\zeta$ small enough that if $x_i(2)/x_i(1) < w + \eta$ and $1 \leq m \leq \varepsilon|z|$, then the cone $\mathcal{C}(B_m, z, \zeta)$ does not intersect $\{y : y(2)/y(1) > \tilde{w} - \eta\}$.

Hence if $\varepsilon$ is small enough and $i$ is large enough, the cone $\mathcal{C}(B_m, x_i, C|x_i|^{-\delta})$ does not intersect with $\{y : y(2)/y(1) > \tilde{w} - \eta\}$.

From Corollary 3.1,

$$R^{\text{out}}[B_m, x_i] \subseteq \mathcal{C}(B_m, x_i, C|x_i|^{-\delta}) \tag{3.20}$$

for all $m \leq \varepsilon|x_i|$ eventually, w.p.1. Thus we indeed have that for $|x_i|, |y_j|$ large,

$$y_j \notin R^{\text{out}}[B_m, x_i]. \tag{3.21}$$

for all $m \leq \varepsilon|x_i|$.

Then for $i$ large enough, for all $m \leq \varepsilon|x_i|$, the set $R^{\text{out}}[B_m, x_i]$ does not intersect with $\{y : y(2)/y(1) > \tilde{w} - \eta\}$.

But if $i$ is large enough, and so $x_i$ is large enough, then, for every $j$ with $|y_j| \geq |x_i|$, we have $y_j \in \{y : y(2)/y(1) > \tilde{w} - \eta\}$. Since $R^{\text{out}}[B_m, x_i]$ contains only points $y$ with $|y| \geq |x_i|$. This completes the proof of Proposition 3.2. $\square$

3.3. *Asymptotic shape.* In this section we prove Proposition 2.1. We are interested in $G(z)$ which is equal to $\max(G_1(z), G_2(z))$, where $G_1(z) = \sup_{1 \leq k \leq z(2)} T(A_k \to z)$ and $G_2(z) = \sup_{1 \leq m \leq z(1)} T(B_m \to z)$. Proposition 2.1 is an immediate consequence of the following proposition, since $p(z)$ defined in (2.18) is the maximum of $p_1(z)$ and $p_2(z)$ defined below.

PROPOSITION 3.4. *Suppose that $\gamma_0$ satisfies* (2.4). *Then with probability* 1,

$$\lim_{|z| \to \infty, z \in \mathbb{Z}_+^2} \frac{G_2(z) - p_2(z)}{|z|} = 0, \tag{3.22}$$

*where*

$$(3.23) \quad p_2(z) = \begin{cases} (\sqrt{z(1)} + \sqrt{z(2)})^2, & \text{if } \frac{z(2)}{z(1)} \geq \left(\frac{\rho}{1-\rho}\right)^2, \\ \frac{z(1)}{1-\rho} + \frac{z(2)}{\rho}, & \text{if } \frac{z(2)}{z(1)} \leq \left(\frac{\rho}{1-\rho}\right)^2, \end{cases}$$

*and similarly*

$$(3.24) \quad \lim_{|z| \to \infty, z \in \mathbb{Z}_+^2} \frac{G_1(z) - p_1(z)}{|z|} = 0,$$

*where*

$$(3.25) \quad p_1(z) = \begin{cases} (\sqrt{z(1)} + \sqrt{z(2)})^2, & \text{if } \frac{z(2)}{z(1)} \leq \left(\frac{\lambda}{1-\lambda}\right)^2, \\ \frac{z(1)}{1-\lambda} + \frac{z(2)}{\lambda}, & \text{if } \frac{z(2)}{z(1)} \geq \left(\frac{\lambda}{1-\lambda}\right)^2. \end{cases}$$

PROOF. We prove (3.22), since (3.24) is the symmetric statement.

We will use the moderate deviations estimate of Lemma 3.1 together with simple deterministic estimates (which we sometimes indicate in outline since their justification is simple but long).

We are interested in $G_2(z) = \sup_{1 \leq m \leq z(1)} T(B_m \to z)$ for $z \in \mathbb{Z}_+^2$.

Recall that $B_m = (m, \beta_m + 1)$, and that from (2.4) we have that

$$(3.26) \quad \frac{\beta_m}{m} \to -\frac{\rho}{1-\rho} := -d_\rho \leq 0.$$

From this convergence, we have that, for any $\tilde{\varepsilon} > 0$, there exists some $K$ such that

$$(3.27) \quad \beta_m + 1 - (-md_\rho) \leq \max(K, \tilde{\varepsilon} m)$$

for all $m \geq 1$. From the form of the function $\mu(z) = (\sqrt{z(1)} + \sqrt{z(2)})^2$, this implies that for any $\varepsilon > 0$,

$$(3.28) \quad |\mu(z - B_m) - \mu(z - (m, -md_\rho))| \leq \varepsilon |z|$$

for all $m \leq z(1)$, for all except finitely many $z \in \mathbb{Z}_+^2$.

Because $\beta_m \leq -1$ for each $m$, and using the convergence in (3.26) again, we have that for some $M$,

$$(3.29) \quad \frac{|z|}{M} \leq |z - B_m| \leq M|z|$$

for all $z \in \mathbb{Z}_+^2$ and all $1 \leq m \leq z(1)$. Then we can apply the moderate deviations estimate of Lemma 3.1 to obtain that for any $\varepsilon > 0$, there exists $\tilde{b} > 0$ such that for all $z \in \mathbb{Z}_+^2$ and $1 \leq m \leq z(1)$,

$$(3.30) \quad \mathbb{P}(|T(B_m \to z) - \mu(z - B_m)| > \varepsilon|z|) \leq \tilde{b} \exp(-\tilde{b}|z|^{1/2}).$$

Summing over $z$ and $m$, and applying Borel–Cantelli, gives the following holds with probability 1: for all except finitely many $z \in \mathbb{Z}_+^2$,

(3.31) $$|T(B_m \to z) - \mu(z - B_m)| \leq \varepsilon|z|$$

for all $1 \leq m \leq z(1)$.

Since $\varepsilon > 0$ is arbitrary, we can combine (3.28) and (3.31) and use $G_2(z) = \sup_{1 \leq m \leq z(1)}$ to give that, with probability 1,

(3.32) $$G_2(z) = \sup_{1 \leq m \leq z(1)} \mu(z - (m, -md_\rho)) + o(|z|)$$

as $|z| \to \infty$ with $z \in \mathbb{Z}_+^2$.

Finally we will estimate the RHS of (3.32). We have

(3.33) $$\sup_{1 \leq m \leq z(1)} \mu(z - (m, -md_\rho))$$
$$= \sup_{1 \leq m \leq z(1)} (\sqrt{z(1) - m} + \sqrt{z(2) + md_\rho})^2$$
$$= z(1) \sup_{x = \frac{1}{z(1)}, \frac{2}{z(1)}, \ldots, 1} \left(\sqrt{1-x} + \sqrt{\frac{z(2)}{z(1)} + xd_\rho}\right)^2,$$

where we put $x = m/z(1)$. As $|z| \to \infty$ with $z \in \mathbb{Z}_+^2$, it can easily be shown that (3.33) is

(3.34) $$z(1) \sup_{0 \leq x \leq 1} \left(\sqrt{1-x} + \sqrt{\frac{z(2)}{z(1)} + xd_\rho}\right)^2 + o(|z|).$$

Calculating this supremum, one finds that it is exactly equal to $p_2(z)$. So combining with (3.32), we have that with probability 1, $G_2(z) = p_2(z) + o(|z|)$, which is what we require. □

### 3.4. Proofs of laws of large numbers.

PROOF OF THEOREM 1 FOR $\lambda > \rho$. Assume $\lambda > \rho$. Let $d_\rho = \rho/(1-\rho)$ as before, and let $d_\lambda = \lambda/(1-\lambda)$, so $d_\rho < d_\lambda$.

LEMMA 3.6. *Let $d_\rho^2 < w < \tilde{w}$. Let $(x_1, x_2, \ldots)$ and $(y_1, y_2, \ldots)$ be infinite directed paths, with $x_i(2)/x_i(1) \to w$ and $y_j(2)/y_j(1) \to \tilde{w}$. Then*

$$\mathbb{P}(\{x_i\} \text{ contains infinitely many points of } \Gamma_\infty^1$$
$$\text{and } \{y_j\} \text{ contains infinitely many points of } \Gamma_\infty^2) = 0.$$

PROOF. Suppose there are infinitely many $j$ such that $y_j \in \Gamma_2^\infty$. For such a $y_j$, its optimal path starts at the point $B_{M(y_j)} = (M(y_j), \beta_{M(y_j)} + 1)$.

Suppose $x_i \in \Gamma_\infty^1$, so that its optimal path starts at some $(\alpha_k + 1, k) = A_k$. Choose $j$ such that $y_j > x_i$ and $y_j \in \Gamma_2^\infty$.

The path from $B_{M(y_j)}$ to $y_j$ cannot pass to the left of $x_i$, since then it would cross the path from $A_k$ to $x_i$. (With probability 1, two maximizing paths which start at different points in $\{A_k, k \geq 1\} \cup \{B_m, m \geq 1\}$ cannot cross; this would contradict the uniqueness of the maximizing path to the point which they share.) Hence the path must pass to the right of $x_i$, and will pass through $x_{i'}$ for some $i' > i$. Then $x_{i'}$ is on the optimal path from $B_{M(y_j)}$ to $y_j$.

If there are infinitely many $i$ with $x_i \in \Gamma_\infty^1$, then there will be infinitely many such $i'$. But from Proposition 3.1, this is an event of probability 0. □

Consider some countable collection of directed paths $(z_i^q, i = 1, 2, \ldots)$, $q \in \mathbb{Q}_+$, each having direction $q$ in the sense that $z_i^q(2)/z_i^q(1) \to q$.

For $q \in \mathbb{Q}_+$ and $r \in \{1, 2\}$, let $P_r^q$ be the property that $\{z_i^q\}$ includes infinitely many points of $\Gamma_r^\infty$.

Let $W = \inf\{q : P_1^q \text{ holds}\}$. If $P_1^q$ holds, then there are infinitely many points in $\{z_i^q\}$ which lie above the competition interface. Let $\tilde{q} > q$; eventually the path $(z_i^q)$ lies above the path $(z_i^{\tilde{q}})$, so $P_1^{\tilde{q}}$ also holds. Thus $P_1^q$ holds for all $q > W$.

LEMMA 3.7.
  (i) If $d_\rho^2 < q < \tilde{q}$, then $\mathbb{P}(P_1^q \text{ and } P_2^{\tilde{q}} \text{ both hold}) = 0$.
  (ii) If $q < \tilde{q} < d_\lambda^2$, then $\mathbb{P}(P_1^q \text{ and } P_2^{\tilde{q}} \text{ both hold}) = 0$.

PROOF. Part (i) follows directly from Lemma 3.6. Part (ii) is the exactly symmetric statement, reversing the roles of $\{A_k\}$ and $\{B_m\}$, of clusters 1 and 2, and of $\rho$ and $\lambda$. □

So suppose $q < \tilde{q}$ with $\tilde{q} - q < d_\lambda^2 - d_\rho^2$. Then either $d_\rho^2 < q < \tilde{q}$ or $q < \tilde{q} < d_\lambda^2$. So using Lemma 3.7 and countable additivity, $\mathbb{P}(P_1^q \text{ and } P_2^{\tilde{q}} \text{ both hold for some such pair } q, \tilde{q}) = 0$.

Hence also $W = \sup\{\tilde{q} : P_2^{\tilde{q}} \text{ holds}\}$, and $P_2^{\tilde{q}}$ holds for all $\tilde{q} < W$.

So for all $q, \tilde{q}$ with $\tilde{q} < W < q$, all but finitely many points in $(z_i^{\tilde{q}})$ lie in $\Gamma_\infty^2$, and all but finitely many points in $(z_i^q)$ lie in $\Gamma_\infty^1$. So the competition interface $\phi_n$ eventually lies between these two paths. Hence indeed $\phi_n$ has direction $W$, and the proof of Theorem 1 in the case $\lambda > \rho$ is complete. □

PROOFS OF THEOREM 1 FOR $\lambda \leq \rho$, AND OF THEOREM 2(b) AND (c).

Suppose $\lambda \leq \rho$. Let $w^* = \frac{\lambda\rho}{(1-\lambda)(1-\rho)}$. Note that $p_1((1, w^*)) = p_2((1, w^*))$; the angle $\alpha$ with $\tan\alpha = w^*$ marks the direction of the shock.

Let $\varepsilon > 0$ and consider the set

$$(3.35) \qquad M_+ = \left\{z \in \mathbb{Z}_+^2 : \frac{z(2)}{z(1)} > \tan(\alpha + \varepsilon)\right\}.$$

If $z \in M_+$, then $p_1(z) > p_2(z)$, and in fact it is easy to show from the form of the functions $p(1)$ and $p(2)$ that this inequality is "uniform" in the following sense: there exists $\delta > 0$ such that for all but finitely many $z \in M_+$,

$$(3.36) \qquad p_1(z) > p_2(z) + \delta|z|.$$

Hence by the limiting shape result Proposition 3.4, with probability 1, for all but finitely many $z \in M_+$ one has $G_1(z) > G_2(z)$, and so $z \in \Gamma_\infty^1$.

Similarly if

$$(3.37) \qquad M_- = \left\{z \in \mathbb{Z}_+^2 : \frac{z(2)}{z(1)} < \tan(\alpha - \varepsilon)\right\},$$

then w.p.1, for all but finitely many $z \in M_-$, $G_2(z) > G_1(z)$ and so $z \in \Gamma_\infty^2$.

Hence the competition interface must eventually lie within the cone of width $2\varepsilon$ around the angle $\alpha$. But $\varepsilon$ is arbitrary, and so Theorem 1 follows for $\lambda \leq \rho$, along with the limiting value given in Theorem 2(a).

An analogous argument implies the result in Theorem 2(b). One just needs to check from the form of $p_1$ and $p_2$ for $\lambda > \rho$ that if

$$(3.38) \qquad \tilde{M}_+ = \left\{z \in \mathbb{Z}_+^2 : \frac{z(2)}{z(1)} > \left(\frac{\lambda}{1-\lambda}\right)^2 + \varepsilon\right\},$$

then for some $\delta > 0$ and all but finitely many $z \in \tilde{M}+$, one has $p_1(z) > p_2(z) + \delta|z|$; and similarly for

$$(3.39) \qquad \tilde{M}_- = \left\{z \in \mathbb{Z}_+^2 : \frac{z(2)}{z(1)} < \left(\frac{\rho}{1-\rho}\right)^2 - \varepsilon\right\},$$

with $p_1$ and $p_2$ reversed. The shape result Proposition 3.4 can then be applied as above. $\square$

## 4. Proofs of fluctuation results.

4.1. *Fluctuations with a concave initial growth interface.* The proof of Theorem 4 follows the same method due to Ferrari and Fontes [6] applied to the vector $\psi(t) = (I(t), J(t))$.

PROOF OF THEOREM 4. The limiting covariances (2.29), (2.30) and (2.31) and the central limit theorem follow from (2.27) and (2.28) and the fact that under $\nu_{\lambda,\rho}$ the variables $N$ are just sum of independent Bernoulli random variables.

Let $F^\zeta_{t,r}$ be the flux of exclusion particles through the space–time line $(0,0)$–$(r,t)$ when the initial configuration is $\zeta$, defined by

(4.1)
$$\begin{aligned} F^\zeta_{t,r} = {}& \text{number of } \zeta \text{ particles to the left of } 0 \text{ at time } 0 \text{ and} \\ & \text{to the right of } r \text{ at time } t \\ & - \text{number of } \zeta \text{ particles to the right of } 0 \text{ at time } 0 \text{ and} \\ & \text{to the left of } r \text{ at time } t. \end{aligned}$$

If $\zeta$ is distributed according to the product measure with density $\bar\rho$ and $r \neq 1 - 2\bar\rho$, Ferrari and Fontes [6] proved that the flux $F^\zeta_{t,r}$ depends on the initial configuration in the following sense:

(4.2)
$$\lim_t \frac{1}{t} \int \nu_{\bar\rho}(d\zeta) \mathbb{E}(F^\zeta_{t,rt} - N^\zeta_{t(r-1+2\bar\rho),0} - \bar\rho^2 t)^2 = 0.$$

Consider a system of first- and second-class particles $(\sigma_t, \xi_t)$ such that $\sigma_t$ is the exclusion process with (marginal) law $\nu_\rho$ and $\sigma_t + \xi_t$ is the exclusion process with marginal law $\nu_\lambda$. Call $\nu_2$ the distribution of $(\sigma_0, \xi_0)$ and $\nu'_2$ this measure conditioned to have a $\xi$ particle at the origin. The $\sigma$ particles are first-class particles and the $\xi$ particles are second-class; see [9]. Let $T_x(\sigma, \xi)(y) = \sigma(y) + \xi(y)\mathbf{1}\{y > x\}$ be the transformation that identifies the $\sigma$ and $\xi$ particles to the right of $x$ and holes and $\xi$ particles to the left of $x$. If $\eta_0 = T_0(\sigma_0, \xi_0)$, then

(4.3)
$$\eta_t = T_{X(t)}(\sigma_t, \xi_t).$$

In [7] it is proved that $X^\eta(t)$, the second-class particle in the exclusion process with initial distribution $\eta$, is the same as $X^{\sigma,\xi}(t)$, a tagged $\xi$ particle with initial configuration $(\sigma, \xi)$. It is also proved that $X^{\sigma,\xi}(t)$ depends on the initial configuration in the following sense:

(4.4)
$$\begin{aligned} \lim_{t \to \infty} \frac{1}{t} \int \nu_2(d(\sigma,\xi)) \\ \times \mathbb{E}(X^{\sigma,\xi}(t) - ((N^\gamma_{0,(\rho-\lambda)t} - N^\sigma_{-(\lambda-\rho)t,0})/(\lambda-\rho)))^2 = 0, \end{aligned}$$

where $\gamma(x) = 1 - \sigma_t(x) + \xi_t(x)$ indicates the holes.

Call $\gamma_t$ the positions of the holes at time $t$: $\gamma_t(x) = 1 - \sigma_t(x) + \xi_t(x)$. Since holes cannot pass from the left to the right of $X(t)$,

(4.5)
$$I(t) = -F^\gamma_{t,X(t)},$$

the negative flux of holes through the line $(0,0)$–$(t, X(t))$. We divide this quantity in two parts:

$$-F^{\gamma}_{t,X(t)} = -F^{\gamma}_{t,vt} + \sum_{x=vt}^{X(t)} \gamma_t(x). \tag{4.6}$$

Apply (4.2) to the hole-process $\gamma_t$ for $r = v$ and $\bar{\rho} = 1 - \rho$ (the density of $\gamma$ particles). Noticing that the holes do exclusion with drift to the left producing a change of sign, we get

$$-F^{\gamma}_{t,vt} \sim N^{\gamma}_{0,(\rho-\lambda)t} + (1-\rho)^2 t \tag{4.7}$$

in the sense of (4.2). On the other hand, by the law of large numbers for the measure as seen from the second-class particle [7], in the same sense,

$$\sum_{x=vt}^{X(t)} \gamma_t(x) \sim (X(t) - vt)(1 - \rho). \tag{4.8}$$

Substituting (4.7), (4.8) and (4.4) in (4.6),

$$\begin{aligned}
-F^{\gamma}_{t,X(t)} &\sim N^{\gamma}_{0,(\rho-\lambda)t} + (1-\rho)^2 t \\
&\quad + \Big(\frac{N^{\gamma}_{0,(\rho-\lambda)t} - N^{\sigma}_{-(\rho-\lambda)t,0}}{\rho - \lambda} - vt\Big)(1-\rho) \\
&= \frac{(1-\lambda)N^{\gamma}_{0,(\rho-\lambda)t} - (1-\rho)N^{\sigma}_{-(\rho-\lambda)t,0}}{\rho - \lambda} + \lambda(1-\rho)t.
\end{aligned} \tag{4.9}$$

Noticing that $N^{\gamma}_{0,(\rho-\lambda)t} = N^{1-\eta}_{0,(\rho-\lambda)t}$ and $N^{\sigma}_{-(\rho-\lambda)t,0} = N^{\eta}_{-(\rho-\lambda)t,0}$ and using (4.5) we get (2.27).

By the same reasons, $J^{\eta}(t) = F^{\sigma}_{t,X(t)}$ which can be computed as before by

$$\begin{aligned}
F^{\sigma}_{t,X(t)} &\sim N^{\sigma}_{-(\rho-\lambda)t,0} + \lambda^2 t - \Big(\frac{N^{\gamma}_{0,(\rho-\lambda)t} - N^{\sigma}_{-(\rho-\lambda)t,0}}{\rho - \lambda} - vt\Big)\lambda \\
&= \frac{-\lambda N^{\gamma}_{0,(\rho-\lambda)t} + \rho N^{\sigma}_{-(\rho-\lambda)t,0}}{\rho - \lambda} + (1-\rho)\lambda t
\end{aligned} \tag{4.10}$$

from where we get (2.28). □

4.2. *Fluctuations in equilibrium.* In this subsection we prove Theorem 5. We consider the stationary exclusion process under the invariant measure $\nu_\rho$, the product measure with density $\rho$. The generator of the process is given by

$$Lf(\eta) = \sum_{x \in \mathbb{Z}} \eta(x)(1 - \eta(x+1))[f(\eta - \delta_x + \delta_{x+1}) - f(\eta)]. \tag{4.11}$$

The measure $\nu_\rho$ is invariant for $L$: $\nu_\rho L = \nu_\rho$. The reverse process with respect to $\nu_\rho$ has generator $L^*$ which is also a totally asymmetric simple exclusion process with reversed jumps:

$$(4.12) \qquad L^* f(\eta) = \sum_{x \in \mathbb{Z}} \eta(x)(1 - \eta(x-1))[f(\eta - \delta_x + \delta_{x-1}) - f(\eta)].$$

A construction of the stationary process with time marginal distribution $\nu_\rho$ can be done by choosing a configuration $\eta$ according to $\nu_\rho$ and then running the process with generator $L$ forward in time and the process with generator $L^*$ backward in time. Let $\underline{\eta} = (\eta_t, t \in \mathbb{R})$ be the stationary process so constructed. The reverse process $\underline{\eta}^*$ is given by $\eta_t^* = \eta_{-t-}$.

The particle jumps of $\underline{\eta}$ induce a stationary point process $S$ in $\mathbb{Z} \times \mathbb{R}$. Let $S_x$ be the (discrete and random) set of times for which a particle of $\underline{\eta}$ jumps from $x$ to $x+1$, and $S = (S_x, x \in \mathbb{Z})$. The map $\underline{\eta} \mapsto S$ associates alternate point processes to each trajectory (see Figure 6, here alternate means that between two successive points in $S_x$ there is exactly one point in $S_{x+1}$). Conversely, for each alternate point configuration $S$ there is a unique trajectory with jump times $S$: $\eta_t(x) = 1$ when the most recent point before $t$ in $S_{x-1} \cup S_x$ belongs to $S_{x-1}$ and 0 otherwise. The law of the process $S$ is space and time translation invariant. Let $S^0$ be the Palm version of $S$, that is, the process with the law of $S$ conditioned to have a point at $(x,t) = (0,0)$. In the corresponding process $\underline{\eta}^0$ there is a particle jumping from 0 to 1 at time zero. In the reverse process $\underline{\eta}^{*0}$ there is a particle jumping from 1 to 0 at time 0.

A *point map* is a function $\pi$ from the point configuration space where the process $S^0$ lives to $\mathbb{Z} \times \mathbb{R}$ such that $\pi(S^0) \in S^0$ with probability 1; see [27]. We shall define a family of point maps. Label the times of $S^0$ by means of a random function $G: \mathbb{Z} \times \mathbb{Z} \to \mathbb{R}$; in fact $G = G(S^0)$, but we drop the dependency of $S^0$ in the notation. Let $G(0,0) = 0$, interpret it as "particle labeled 0 jumps to hole labeled 0 at time 0." This determines the label of all other times in $S^0$ as follows: the time $G(i,i)$ is the $i$th positive time in $S_0^0$ for positive $i$'s and the $(-i)$th negative time for negative $i$'s:

$$(4.13) \quad \begin{aligned} G(i,i) &:= \inf(S_0^0 \cap (G(i-1,i-1), \infty)) \\ &= \sup(S_0^0 \cap (-\infty, G(i+1,i+1))). \end{aligned}$$

Assuming $G(i,j)$ and $G(i+1,j+1)$ are given, define

$$(4.14) \qquad G(i+1,j) = S^0_{i+1-j} \cap (G(i,j), G(i+1,j+1))$$

(which almost surely is a unique point) and

$$(4.15) \qquad G(i,j+1) = S^0_{i-j-1} \cap (G(i,j), G(i+1,j+1))$$

(which is also a unique point almost surely).

To interpret $G$ in function of the particle motion, label the particles of $\eta_0^0$ in decreasing order, giving the label 0 to the particle at site 1. Call $P_j(0)$ the position of the $j$th particle at time zero; we have $P_0(0) = 1$ and $P_{j+1}(0) < P_j(0)$, $j \in \mathbb{Z}$. Label the holes of $\eta_0^0$ in increasing order, giving the label 0 to the hole at site 0: $H_0(0) = 0$ and $H_{i+1}(0) > H_i(0)$ for all $i$. The position of the $j$th particle and the $i$th hole at time $t$ are, respectively, $P_j(t)$ and $H_i(t)$. The order is preserved at later and earlier times: $P_j(t) > P_{j+1}(t)$ and $H_i(t) < H_{i+1}(t)$, for all $t \in \mathbb{R}$, $i,j \in \mathbb{Z}$. At time $G(i,j)$ the $i$th hole and the $j$th particle of $\eta^0$ interchange positions; in particular $G(0,0) = 0$. Let $\underline{G} = (G(z), z \in \mathbb{Z}^2)$. Since $\underline{G}$ is a deterministic function of $\eta^0$ we write when necessary $\underline{G}(\eta^0)$.

The random function $G$ induces a family of point maps $\pi_{i,j}$ given by

$$(4.16) \qquad \pi_{i,j}(S^0) = (i - j, G(i, j)).$$

The space–time shift by $(z,t) \in \mathbb{Z} \times \mathbb{R}$ of the point process $S$ is defined by $S - (z,t) = \{(z' - z, t' - t), (z', t') \in S\}$. Each map $\pi_{i,j}$ is bijective in the sense that

$$(4.17) \qquad (z,t) \mapsto \pi_{i,j}(S^0 - (z,t)) + (z,t)$$

is a bijection in $S^0$. Or, said in other words, each point of $S^0$ is mapped to another unique point in $S^0$ by means of $\pi$; distinct points are mapped to distinct points.

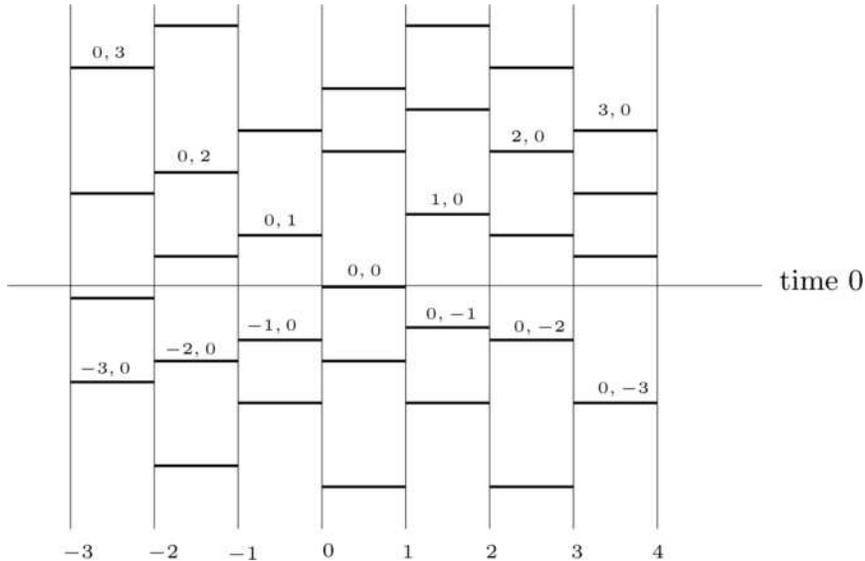

Fig. 6. Labeling of $S^0$. The horizontal axis is $\mathbb{Z}$; time runs upwards. The heights of the horizontal segments between vertical lines $x$ and $x+1$ represent the times $G(i,j)$ with $i - j = x$.

LEMMA 4.1. *The law of*

(4.18) $$S^0 - (\bar{i} - \bar{j}, G(\bar{i}, \bar{j}))$$

*does not depend on* $\bar{z} = (\bar{i}, \bar{j}) \in \mathbb{Z}^2$. *Consequently, the law of*

(4.19) $$[G(z - \bar{z}) - G(\bar{z}), z \in \mathbb{Z}^2]$$

*does not depend on* $\bar{z}$.

PROOF. The second statement is an immediate consequence of the first one.

Let $\Theta_{i,j}$ be the map from the alternate point-configuration space onto itself defined by

(4.20) $$\Theta_{i,j}(S^0) = S^0 - \pi_{i,j}(S^0).$$

$\Theta_{i,j}$ maps the point-configuration $S^0$ to the same configuration shifted by (or "as seen from") $\pi_{i,j}(S^0) = (i - j, G(i, j))$. This map is invertible, its inverse is $\Theta_{-i,-j}$. Then Theorem 3.1 of Heveling and Last [14] implies $S^0$ and $S^0 - (i - j, G(i, j))$ are identically distributed for each $(i, j)$. □

Define the families $\underline{X} = (X(i,j), i,j \in \mathbb{Z})$ and $\underline{Y} = (Y(i,j), i,j \in \mathbb{Z})$ by

(4.21) $$X(i,j) = G(i,j) - \max(G(i, j-1), G(i-1, j)),$$

(4.22) $$Y(i,j) = \min(G(i, j+1), G(i+1, j)) - G(i,j).$$

The particle interpretation of $G(i,j)$ says that $X(i,j)$ is the time the particle $j$ waits to jump over hole $i$ starting at the time they become neighbors (i.e., the maximum of the time particle $j - 1$ jumped over hole $i$ and the time particle $j$ jumped over hole $i - 1$). Since $-G(i,j)$ are the jumping times of particles for the reverse process, which has the same law as the direct process but with reversed jumps, $Y(i,j)$ has the same interpretation as $X(i,j)$ for the reverse process and hence the same law.

LEMMA 4.2. *The variables* $(X(i,j), i,j \in \mathbb{Z})$ *are independent and identically distributed with exponential law of mean* 1. *The same is true for* $(Y(i,j), i,j \in \mathbb{Z})$.

PROOF. First take the set of labels $(i,j)$ such that $j > 1$ and $i > 1$. The pairs $(i,j)$ in this set satisfy $G(i-1, j) > 0$ and $G(i, j-1) > 0$ with probability 1. The conditioning event $\{G(0,0) = 0\}$ is measurable for the sigma field generated by $(\eta_s, s \leq 0)$. Then the Markov property implies that the waiting times $X(i,j)$ with $(i,j)$ in the above set are independent and exponentially distributed of parameter 1.

To extend the result to all $(i,j)$, consider a finite set of indexes $I \subset \mathbb{Z}^2$. It suffices to show that $(X(z), z \in I)$ are i.i.d. exponential of parameter 1. Take $\bar{z} = (\bar{i}, \bar{j})$ such that $j - 1 > \bar{j}$ and $i - 1 > \bar{i}$ for all $(i,j) \in I$. Then, with probability 1 $G(\bar{z}) < \min\{G(i-1,j), G(i,j-1)\}$ for all $(i,j) \in I$. Apply the same reasoning as above for the family $(G(z - \bar{z}) - G(\bar{z}), z \in \mathbb{Z}^2)$ which has the same law as $\underline{G}$ by Lemma 4.1. $\square$

For each $z \in \mathbb{Z}^2$ let

$$(4.23) \qquad z^{\text{NE}}(z, \underline{G}) := \arg\min\{G(z + (0,1)), G(z + (1,0))\}$$

be the point in $\{z + (0,1), z + (1,0)\}$ that realizes the minimum between $G(z + (0,1))$ and $G(z + (1,0))$ [recall (2.12)]. A path $(z_n, n \geq 0)$ such that $z_{n+1} = z^{\text{NE}}(z_n, \underline{G})$ is called a *North-East Minimizing Sequence* (NEMS) for $\underline{G}$.

LEMMA 4.3. *Suppose $\underline{G}$ and $\underline{Y}$ satisfy (4.22). If a path $(z_n, n \geq 0)$ is a NEMS for $\underline{G}$, then it is a geodesic for $\underline{Y}$.*

PROOF. First consider $\underline{X}$ related to $\underline{G}$ by (4.21). For each $z \in \mathbb{Z}^2$, let

$$(4.24) \qquad z^{\text{SW}}(z, \underline{G}) := \arg\max\{G(z - (0,1)), G(z - (1,0))\}$$

be the point in $\{z - (0,1), z - (1,0)\}$ that realizes the maximum between $G(z - (0,1))$ and $G(z - (1,0))$. Let $z_0, z_1, \ldots, z_m$ be an up-right path. If $z_l = z^{\text{SW}}(z_{l+1}, \underline{G})$ for all $l = 0, 1, \ldots, m-1$, then it follows that the path $z_0, z_1, \ldots, z_m$ is the geodesic between $z_0$ and $z_m$, in terms of the variables $\underline{X}$.

Now define for each $i, j$

$$\tilde{G}(i,j) = -G(-i,-j)$$

and

$$\tilde{X}(i,j) = \tilde{G}(i,j) - \max\{\tilde{G}(i,j-1), \tilde{G}(i-1,j)\}.$$

Since (4.23) holds, it follows that $-z_{l+1} = z^{\text{SW}}(-z_l, \tilde{\underline{G}})$ for any $l$, so that by the previous paragraph, the path $-z_m, -z_{m-1}, \ldots, -z_0$ is the geodesic between $-z_m$ and $-z_0$ in terms of the weights $\tilde{\underline{X}}$.

Finally note that from (4.22) we get $Y(i,j) = \tilde{X}(-i,-j)$ for all $i, j$. Thus $z_0, \ldots, z_{m-1}, z_m$ is the geodesic between $z_0$ and $z_m$ in terms of the weights $\underline{Y}$. Thus by definition the path $(z_n, n \geq 0)$ is a semi-infinite geodesic for the weights $\underline{Y}$, as desired. $\square$

The following result is proved in [20]:

LEMMA 4.4. *Let $\underline{Y} = (Y(i,j), i,j \in \mathbb{Z})$ be i.i.d. exponential with mean 1. Then for any $\alpha \in [0, \pi/2]$, there is a.s. only one geodesic for the weights $\underline{Y}$ starting at the origin with asymptotic direction $\alpha$.*

PROOF OF THEOREM 5. From (2.12), the competition interface is a NEMS for $\underline{G}$. Hence by Lemma 4.3, it is a geodesic for $\underline{Y}$, and from Theorem 2 it has direction $\alpha$ where $\tan(\alpha) = \rho^2/(1-\rho)^2$. Hence the competition interface is a.s. equal to the unique geodesic for $\underline{Y}$ in this direction given by Lemma 4.4, and so the fluctuation exponent of the competition interface is the same as that of the geodesic in direction $\alpha$, so that $\chi = \xi_1(\alpha)$ as desired. □

4.3. *Fluctuations with a convex initial growth interface.* From the competition interface $\phi = (\phi_n)_{n \geq 0}$, we can obtain two subsequences $(\phi_{n_l^A})_{l>0}$ and $(\phi_{n_l^B})_{l>0}$ such that the optimal path from $\Gamma_0$ to $\phi_{n_l^A}$ comes from $(A_k)_{k>0}$ and the optimal path from $\Gamma_0$ to $\phi_{n_l^B}$ comes from $(B_m)_{m>0}$. Given $z \notin \Gamma_0$, we denote by $A(z)$ the point in $(A_k)_{k>0}$ from which the optimal path from $(A_k)_{k>0}$ comes. Analogously, we denote by $B(z)$ the point in $(B_m)_{m>0}$ from which the optimal path from $(B_m)_{m>0}$ comes.

LEMMA 4.5. $\mathbb{P}$-*a.s., there exist two semi-infinite geodesics* $\pi_A \subseteq \Gamma_\infty^1$ *and* $\pi_B \subseteq \Gamma_\infty^2$ *such that*

$$(4.25) \quad \lim_{l \to \infty} \pi(A(\phi_{n_l^A}), \phi_{n_l^A}) = \pi_A \quad \text{and} \quad \lim_{l \to \infty} \pi(B(\phi_{n_l^B}), \phi_{n_l^B}) = \pi_B.$$

*Furthermore, both semi-infinite geodesics have the same asymptotic direction as the competition interface.*

PROOF. Since $A(\phi_{n_{l+1}^A})$ is always to the northeast of $A(\phi_{n_l^A})$, $(A(\phi_{n_l^A}))_{l \geq 1}$ is a monotone sequence, and so it converges. By Proposition 3.1 the limit must be a point $z_A$, that is,

$$(4.26) \qquad A(\phi_{n_l^A}) = z_A \qquad \text{for all large } l.$$

Define the "rightmost" semi-infinite geodesic $\pi_A := (y_k)_{k>0}$ by the following rule: set $y_1 := z_A$; given $y_k$, then if $y_k + (0,1)$ belongs to some geodesic connecting $z_A$ to some point in $\phi$ and $y_k + (1,0)$ not, then set $y_{k+1} := y_k + (0,1)$, otherwise set $y_{k+1} := y_k + (1,0)$. Notice that there is no semi-infinite geodesic caught between $(y_k)_{k>0}$ and $\phi$, and together with (4.26) and Proposition 7 of Ferrari and Pimentel [11], this implies that $\pi(A(\phi_{n_l^A}), \phi_{n_l^A})$ must converge to $\pi_A$ and have the same asymptotic direction $e^{i\theta}$ of $\phi$. The proof of the convergence for $\pi(B(\phi_{n_l^B}), \phi_{n_l^B})$ follows the same argument. □

PROOF OF THEOREM 6. Denote by $z_A$ the starting point of the semi-infinite geodesic $\pi_A$ and by $z_B$ the starting point of the semi-infinite geodesic $\pi_B$. Since the competition interface is caught in between $\pi_A$ and $\pi_B$,

$$(4.27) \quad |y_{\phi_0}(\theta, r) - re^{i\theta}| \leq \max\{|y_{z_A}(\theta, r) - re^{i\theta}|, |y_{z_B}(\theta, r) - re^{i\theta}|\},$$

although

$$(4.28) \qquad |y_{z_A}(\theta, r) - re^{i\theta}| \leq |y_{z_A}(\theta, r) - (z_A + re^{i\theta})| + |z_A|$$

and

$$(4.29) \qquad |y_{z_B}(\theta, r) - re^{i\theta}| \leq |y_{z_B}(\theta, r) - (z_B + re^{i\theta})| + |z_B|.$$

Together with (4.27) this yields that $\chi \leq \xi_2$. $\square$

**Acknowledgments.** We thank Hermann Thorisson for useful advice on point stationarity and Timo Seppäläinen for nice discussions. We thank the referees for their very careful reading of the paper and their helpful comments and corrections.

P. A. FERRARI
INSTITUTO DE MATEMÁTICA E ESTATÍSTICA
UNIVERSIDADE DE SÃO PAULO
CAIXA POSTAL 66281
05311-970, SÃO PAULO
BRASIL
E-MAIL: [pablo@ime.usp.br](pablo@ime.usp.br)

J. B. MARTIN
DEPARTMENT OF STATISTICS
UNIVERSITY OF OXFORD
1 SOUTH PARKS ROAD
OXFORD OX1 3TG
UNITED KINGDOM
E-MAIL: [martin@stats.ox.ac.uk](martin@stats.ox.ac.uk)

L. P. R. PIMENTEL
DEPARTMENT OF APPLIED MATHEMATICS
  (DIAM)
DELFT UNIVERSITY OF TECHNOLOGY
MEKELWEG 4
2628 CD DELFT
THE NETHERLANDS
E-MAIL: [L.PintoRodriguesPimentel@tudelft.nl](L.PintoRodriguesPimentel@tudelft.nl)